\documentclass[3p,12pt,authoryear]{elsarticle}

\usepackage{amssymb}
\usepackage{mathtools}

\usepackage{hyperref}
\usepackage{float}
\usepackage[misc]{ifsym}
\usepackage{subcaption}
\usepackage{graphicx}
\usepackage{multirow}
\usepackage{longtable}
\usepackage{algorithm2e}
\usepackage[table]{xcolor}
\usepackage{xcolor}
\newcolumntype{P}[1]{>{\centering\arraybackslash}p{#1}}

\journal{Computers \& Chemical Engineering}

\begin{document}

\begin{frontmatter}



\title{Stable optimisation-based scenario generation via game theoretic approach}


\author[1]{Georgios L. Bounitsis}

\author[1]{Lazaros G. Papageorgiou}

\author[1]{Vassilis M. Charitopoulos \corref{cor1}}
\ead{v.charitopoulos@ucl.ac.uk}

\affiliation[1]{organization={The Sargent Centre for Process Systems Engineering, Department of Chemical Engineering, UCL (University College London)},
            addressline={Torrington Place}, 
            city={London},
            postcode={WC1E 7JE}, 
            country={UK}}

\cortext[cor1]{Corresponding author}


\begin{abstract}
Systematic scenario generation (SG) methods have emerged as an invaluable tool to handle uncertainty towards the efficient solution of stochastic programming (SP) problems. The quality of SG methods depends on their consistency to generate scenario sets which guarantee stability on solving SPs and lead to stochastic solutions of good quality. In this context, we delve into the optimisation-based Distribution and Moment Matching Problem (DMP) for scenario generation and propose a game-theoretic approach which is formulated as a Mixed-Integer Linear Programming (MILP) model. Nash bargaining approach is employed and the terms of the objective function regarding the statistical matching of the DMP are considered as players.   Results from a capacity planning case study highlight the quality of the stochastic solutions obtained using MILP DMP models for scenario generation. Furthermore, the proposed game-theoretic extension of DMP enhances in-sample and out-of-sample stability with respect to the challenging problem of user-defined parameters variability.

\end{abstract}


\end{frontmatter}


\section{Introduction}
\label{sec:introduction}

Optimisation under uncertainty and the developments on the corresponding mathematical frameworks constitute a topical domain of the Process Systems Engineering (PSE) literature \citep{Li2021}. Various mathematical approaches can be exploited depending on the characterisation of uncertainty and the degree of risk aversion of the problem at hand. Stochastic Programming (SP) is a risk-neutral approach which exploits scenario-based formulations in order to optimise the expected value of the problem over a known probability distribution. Thus, scenario generation (SG) or scenario reduction (SR) approaches have attracted particular interest. By definition these methods aim to create a smaller and representative set of scenarios to efficiently solve computationally challenging stochastic programs. A wide variety of methods such as copula sampling, machine learning and optimisation models can be exploited for the development of SG frameworks \citep{Kaut2011,Medina-Gonzalez2020,Li2014}. In particular, optimisation-based techniques \citep{Lohndorf2016,Bertsimas2022} constitute a main category of SG methods and include the well-known Moment Matching Problem (MMP) \citep{Høyland2001,Høyland2003}. MMP generates scenario sets by solving a statistical errors’ minimisation problem, which inherently translates into as nonlinear (and nonconvex) programming (NLP) problem. \citet{Calfa2014} proposed an enhancement on the MMP NLP problem by matching in parallel the cumulative probability distributions of the considered parameters. This problem is mentioned as Distribution and Moment Matching Problem (DMP). Although DMP may improve performance and statistical matching of the MMP, it remains nonlinear and nonconvex and may suffer from under-specification issues which worsen its performance. \citet{Bounitsis2022} reformulated DMP as a Mixed-Integer Linear Programming (MILP) model for scenario reduction and proposed its integration in a data-driven scenario generation framework including copula sampling and clustering. This SG methodology was shown to overcome the so-called under-specification issues of the NLP based counterparts. However, results indicated sensitivity of MMP and DMP models against the user-defined parameters regarding the errors' weights of the objective functions. In other words, the solution of MMP and DMP optimisation problem is not robust against the user-defined errors' weights undermining their stability and efficacy.

In this work, an extension to the work by \citet{Bounitsis2022} is proposed aiming to mitigate the impact of the model’s user-defined parameters on the DMP MILP scenario reduction model and consequently to enhance stability and performance of the SG framework. In particular, the DMP MILP model is modified and a Nash bargaining approach is used for the terms of the objective function, which are considered as the players of a game. Game theoretic approaches have been widely employed for applications of the Process Systems Engineering literature and various frameworks to efficiently handle the computational complexity of game theoretic problems are proposed \citep{Faisca2009,Marousi2023}. Thus, DMP problem at hand is modelled as MILP following a separable programming reformulation of the game theoretic problem for the approximation of the Nash product \citep{Gjerdrum2001,Charitopoulos2020}. Finally, the enhanced quality and stability of the optimisation-based MILP DMP model against different sets of user-defined errors' weights are validated through bias, in-sample and out-of-sample tests \citep{Kaut2007}. 

The remainder of the article is organised as follows: in Section \ref{sec:literature} a summary of the main theoretical aspects is provided while the detailed methodology is outlined in Section \ref{sec:methodology}. In Section \ref{sec:casestudy}, the proposed framework is employed to evaluate its stability on a capacity planning case study. Finally, conclusions are drawn in Section \ref{sec:conclusions}.

\section{Preliminaries \& Literature review}
\label{sec:literature}

\subsection{Stochastic Programming}
\label{sec:sp}

Stochastic Programming (SP) constitutes a well-established mathematical framework for optimisation under uncertainty \citep{Sahinidis2004,Li2021}. On its traditional form it is a risk neutral approach and the uncertainty is modelled via a known discrete probability distribution. Compared to robust optimisation,  which optimises over the worst-case scenario \citep{Ben-Tal2002}, stochastic programming leads to solutions that are optimal considering the whole uncertain set, while an expected value is optimised \citep{Birge2011,King2012}. Chance-constrained programming, in opposition with stochastic programming, is risk-averse and uses prespecified measures to quantify the risk level \citep{Li2008}. Moreover, in multi-parametric programming an optimisation model is solved for a range and as a function of multiple uncertain parameters \citep{Oberdieck2016, Charitopoulos2018}. Although different frameworks are specialised for different types of problems, the quantification of uncertainty remains a crucial problem for various frameworks.

Two-stage stochastic programming (TSSP) constitutes the typical version of a stochastic problem. In this approach uncertainty is considered to be realised at one time step in the future. The two stages of the TSSP involve the certain first stage (before the realisation of uncertainty) and the uncertain second stage in the future. Thus, two types of variables are defined for TSSP: (i) first-stage (or "here-and-now") decisions which are determined before the realisation of the uncertain parameters and are independent of the uncertainty, and (ii) the second-stage (or "wait-and-see", or recourse actions) decisions which depend on the realised uncertain parameters, induce corrective impact to the decisions and can alleviate any arising infeasibility. TSSP aims to optimise the objective function of the first-stage costs while optimising the expected value of the second-stage costs \citep{Birge2011}. A mathematical representation of a linear two-stage stochastic programming problem is provided by eqs. \eqref{eq:tssp} - \eqref{eq:expcost} \citep{Shapiro2014}:

\begin{align} 
\label{eq:tssp}
\begin{split}
 \underset{{x \in \mathbb{R}^{N}_{+}}} {\mathrm{min}} \qquad & c^\top \cdot x + \mathbb{E} \left[ Q \left(x,\mathbb{\xi}\right) \right]  \\
& subject \; to \\
& A \cdot x \le b \\
& x \geq 0
\end{split}
\end{align}

where $Q(x,\xi)$ is the optimal value of the second-stage problem:

\begin{align}
\label{eq:expcost}
\begin{split}
 \underset{{y \in \mathbb{R}^{M}_{+}}} {\mathrm{min}} \qquad & q^\top \cdot y \\
& subject \; to \\
& T \cdot x + W \cdot y \le h \\
& y \geq 0
\end{split}
\end{align}

In this formulation $x \in \mathbb{{R}^{N}_{+}}$ denotes the first-stage decisions, $y \in \mathbb{{R}^{M}_{+}}$ is the set of the second-stage decisions, and the vector $\mathbb{\xi} = (q,T,W,h)$ denotes the (known or uncertain) data of the second stage problem. By definition TSSP aims to optimise the expected value over a probability distribution. The distribution of the uncertain parameters, $\mathbb{\xi}$, of problem in Eqs. \eqref{eq:tssp} - \eqref{eq:expcost} can be either discrete or continuous. When the distribution $\mathbb{\xi}$ is discrete, then a finite number of realisations for the uncertain parameters can be considered to describe the probability distribution. In other words, $\mathbb{\xi}$ has a finite support. If $\mathbb{\xi}$ has a finite number of $K$ possible realisations, these realisations are also called \textit{scenarios}, say ${\xi}_k=(q_k,T_k,W_k,h_k)$, with respective probabilities of occurrence $p_k$, $k \in \{1,\dots,K\}$ (scenario set $K=\{1,\dots,K\}$). In such way the standard formulation for the TSSP arises, which can render it computationally tractable. The deterministic equivalent problem for the linear TSSP of Eqs. \eqref{eq:tssp} - \eqref{eq:expcost} is given by Eq. \eqref{eq:deteq}:

\begin{align}\label{eq:deteq}
\begin{split}
 \underset{{x,y_1,...,y_K}} {\mathrm{min}} \qquad & c^\top \cdot x + \sum_{k \in K} p_k \cdot q^\top_k \cdot y_k \\
& subject \; to \\
& A \cdot x \le b \\
& T_k \cdot x + W_k \cdot y_k \le h_k \qquad \forall k \in K \\
& x,\ y_k \geq 0 \qquad \forall k \in K 
\end{split}
\end{align}

\begin{figure}[H]  
\centering
\includegraphics[width=12.4cm]{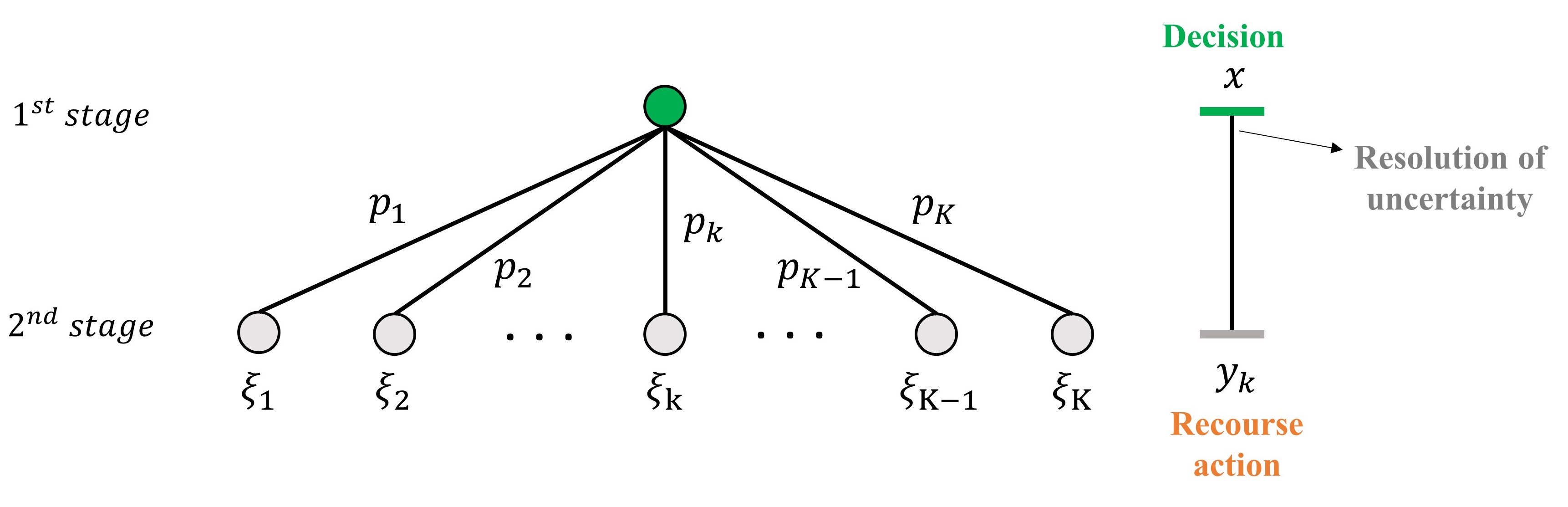}
\caption{Two-stage stochastic problem and scenario set representation, where $x$ are the first-stage decisions, $y_k$ are the second-stage decisions for each scenario $k \in K$, with values $\xi_k$ and probabilities $p_k$.}
\label{fig:sg-tssp}
\end{figure}

In the formulation of Eq. \eqref{eq:deteq}, every scenario $\xi_k=(q_k,T_k,W_k,h_k),\; k \in K =\{ 1, ..., K \}$, results to a two-stage decision vector $y_k$ and by solving the two-stage problem the optimal first-stage decisions $x$ can be computed. In other words, given $x$, each $y_k$ gives the corresponding optimal second-stage decisions for the realisation $\xi_k$, i.e., in scenario $k$. Regarding the nomenclature, matrix $W_k$ is referred to as recourse matrix and matrix $T_k$ is referred to as technology matrix in the literature. As it is demonstrated by Eq. \eqref{eq:deteq} uncertainties can be considered on both of the latter as well as on the right-hand side (RHS) of the constraints or the coefficients of the objective function. A visualisation of the stochastic process and the scenario set is presented in Figure \ref{fig:sg-tssp}.

\subsection{Scenario Generation for Stochastic Programming}
\label{sec:sg}

Based on Eq. \eqref{eq:deteq} is clear that stochastic programs can be computationally challenging or intractable due to either the nature of the optimisation problem (e.g., large number of binary variables in MILP problems) or the true stochastic process. Applications in SP deal with stochastic parameters which are approximated through discrete probability distributions. In case that such discrete distributions are not available or there is a very large set of original data, then available data are processed towards the creation of a smaller and representative subset of scenarios to efficiently solve the stochastic programs. Scenario Generation methods aim to create a set of finite realisations for the uncertain parameters, defining their values and probabilities of occurrence, which is representative of the original stochastic process \citep{King2012}. Nonetheless, methodologies which reduce the size of the scenario set by selecting scenarios and redefining their probabilities are referred to as Scenario Reduction methods. In Fig. \ref{fig:sg-tssp} an instance of a generated scenario set of $K$ scenarios with values $\xi_k$ and probabilities $p_k$ can be envisaged.

\citet{Li2020} offers a literature review of scenario generation methods specialised for wind power data. According to latter review, scenario generation techniques can be classified into three main classes: (i) sampling-based, (ii) forecasting-based, and (iii) optimisation-based. Adopting this classification, the current work delves into the study of the optimisation-based SG techniques. These methodologies are of particular interest in the Operations Research and the Process Systems Engineering literature and has been showcased in critical recent research studies such as the empirical analysis of scenario generation methods by \citet{Lohndorf2016} and the recently proposed optimisation-based scenario reduction method by \citet{Bertsimas2022}.  Two main branches of optimisation-based techniques regard the distance matching problems and the moment matching problem (MMP). An overview of basic and recent methodologies and their variants is presented as follows.

\subsection{Moment Matching Problem (MMP) for Scenario Generation}
\label{sec:mmp}

Generation of two-stage scenario sets and multi-stage scenario trees can be based on solving the Moment Matching Problem (MMP), which was originally introduced by \citet{Høyland2001}. A follow-up work proposing a heuristic algorithm for moment matching problem was presented by \citet{Høyland2003}. The purpose of the MMP given a structure of the scenario tree, i.e., the number of the nodes in every stage, lies in the determination of the optimal values of the uncertain parameters in each node as well as their corresponding probabilities of realisation. MMP achieves that through an error minimisation problem, which generates scenarios by minimising the errors between an original distribution (estimated beforehand) and the reduced set regarding statistical moments of the parameters and their corresponding correlation matrix.

Originally, MMP is formulated as a nonlinear and nonconvex programming (NLP) optimisation problem. When the Euclidean distance ($L^2$-norm) is utilised to quantify the errors then the problem converts into a squared error minimisation problem. However, various reformulations of the MMP can be employed using the Manhattan distance ($L^1$-norm) or the Chebyshev distance ($L^{\infty}$-norm). These reformulations are beneficial as can render the initial NLP problem to a less complicated LP or MILP problem in which the values of the nodes are estimated a priori through a sampling procedure or a clustering method \citep{Xu2012}. Then, the only decision variables are the selection of nodes and/or their probabilities. In case of $L^2$-norm-based formulation of the objective function the Moment Matching Problem can be written as follows:

\begin{align}\label{eq:mmp}
\begin{split}
\underset{x,p \in \mathbb{{R}^{K}_{+}}} {\mathrm{min}} \qquad & \sum_{m\in M}  w_{m} \cdot \left[ f_{m}\left(x,p\right)-{Sval}_{m} \right] ^2
\\ & subject \ to
\\ & \sum_{k=1}^{K} p_{k} = 1
\\ & p_k\ \in\ \left[0,1\right] \qquad \forall k \in K=\{1,...,K\}
\end{split}
\end{align}
where, $x$ is a vector of the values of uncertain parameters, $p$ is a vector of probabilities of nodes, $m \in M$  is the set of statistical properties to be matched, $f_m\left(x,p\right)$ is the mathematical expression of statistical property $m$ calculated from the generated scenarios $k \in K$ with values $x$ and corresponding probability $p$ and ${Sval}_m$ is the value of statistical property $m$ as estimated beforehand by the data/distribution of each parameters \citep{Gulpinar2004}.

\begin{figure}[H] %
    \centering
    \subfloat[\centering Same values may be attributed to several of the prespecified nodes.] 
    {\label{fig:us1} \includegraphics[width=7.5cm]{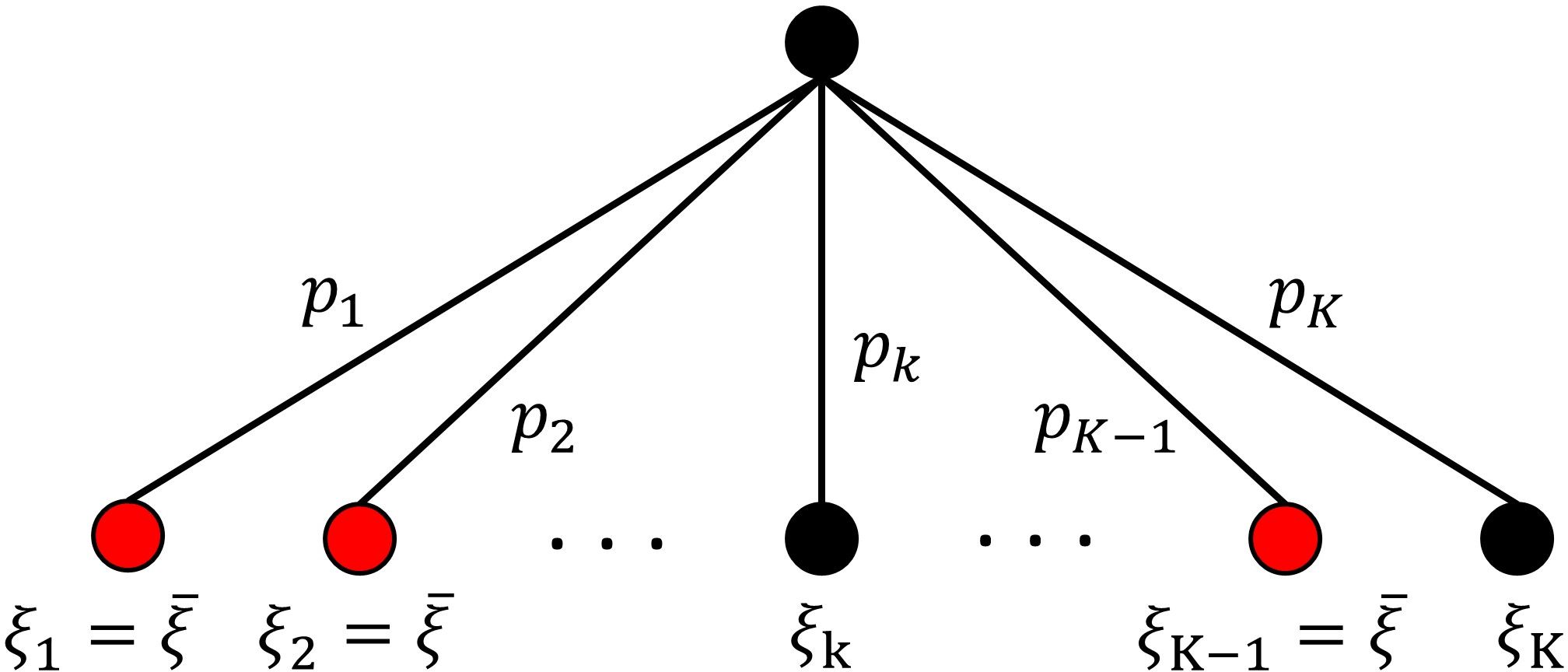} }%
    \qquad
    \subfloat[\centering Zero probability of occurrence may be defined for several of the prespecified nodes.]
    {\label{fig:us2} \includegraphics[width=7.5cm]{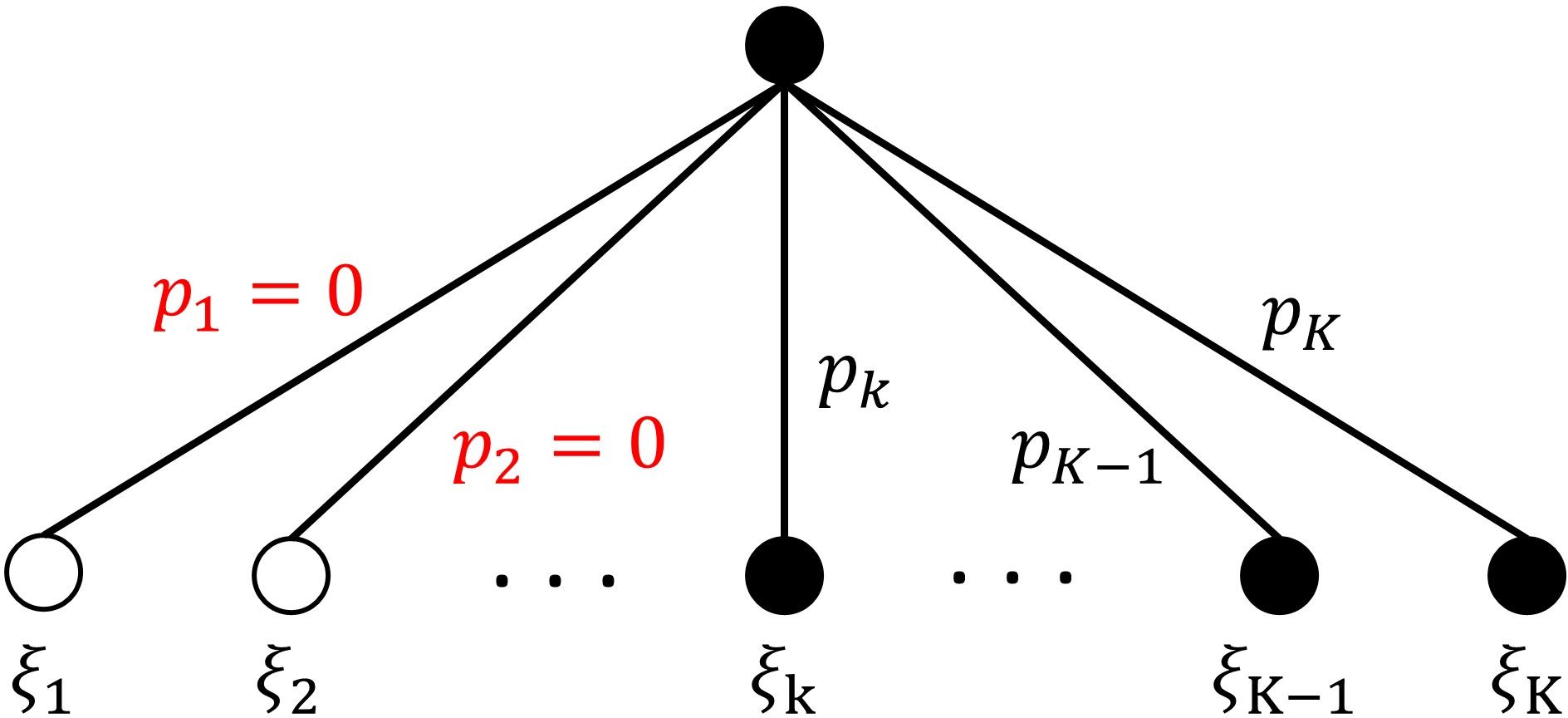} }%
    \caption{Types of under-specification issues that may arise from the solution of original NLP MMP problems and worsen their performance.}%
    \label{fig:us}%
\end{figure}

\citet{Calfa2014} introduced the Distribution \& Moment Matching Problem (DMP) which aims to parallel match the stochastic distribution of the uncertain parameters, by minimising the errors regarding the empirical cumulative distribution function (ECDF) between original and final distributions. The proposed models remain NLP and have enhanced performance compared to MMP. However, these may lead to under-specification issues, in which case either a unique scenario is assigned to several nodes or zero probabilities are assigned to some nodes, and ultimately the performance of the model may be worsened. A visualisation of possible under-specification issues is presented in Fig. \ref{fig:us}. Recently, \citet{Kaut2021} proposed an alternative approach for the solution of MMP as scenario reduction problem by formulating MMP as MILP and using binomial expansion in order to compute the moments of the generated reduced sets. This formulation uses the known values of the nodes and employ binary variables to indicate the selection of the nodes for the final reduced set. \citet{Bounitsis2022} reformulated DMP as MILP using $L^1$ and $L^{\infty}$ norms. The scenario reduction model was integrated into a data-driven framework with copula-based sampling of original scenarios and clustering techniques to reduce computational complexity on the MILP model. Its unified impact seems to lead to significant mitigation of the under-specification issues of the NLP model and enhanced quality of stochastic solutions. 

Beyond the aforementioned optimisation models, instances of MMP integration to iterative algorithms towards scenario generation have been proposed. \citet{Ji2005} were the first to propose an LP MMP problem towards efficient multi-stage scenario generation. In particular, at each stage the values of the nodes are selected a priori and the LP model redefines their probabilities. \citet{Xu2012} also presented a framework for multi-stage scenario tree generation integrating sequential simulations, K-means clustering, time series and solution of LP MMP highlighting efficiency in capturing inter-stage dependencies and the time-varying fluctuations. \citet{Chopra2020} incorporated dimensionality reduction of original dataset a priori, using principal component analysis, and so the final MMP model resulted to reduced computation times. Besides, decomposition techniques can be combined with MMP improving its computational performance. \citet{Li2016-CHO} presented the integration of Cholesky decomposition and clustering to MMP, while \citet{Mehrotra2013} proposed an optimisation-based decomposition algorithm using cubature formulas and column generation to generate moment matching scenarios.

\subsection{Scenario reduction via probabilistic distance minimisation}
\label{sec:sr}

Scenario reduction methods typically select some representative scenarios of an original uncertain set and determine their probabilities through probabilistic distance minimisation problems. To this end, the minimisation of a distance metric is set as the objective function. For instance, the Wasserstein distance constitutes a well-known distance used for optimal discretisation \citet{Pflug2001}. Work by \citet{Dupacova2003}, which proposed two SR algorithms: (1) forward selection and (2) backward selection, seeks to minimise the global probabilistic distance using canonical probability metric for each scenario individually. \citet{Heitsch2003,Heitsch2007,Heitsch2009} through a series of studies enhanced the notion of the latter algorithms by considering the whole set of the original scenarios at each iteration of their algorithm and using Fortet-Mourier metrics as distance. Overall, they lead to an improved computational performance of the scenario reduction algorithms.

In case that the $L^2$-norm-based Wasserstein distance between two discrete distributions is considered, a linear representation of the mass transportation problem towards scenario reduction can be written as follows:

\begin{align}\label{eq:wm}
\begin{split}
\underset{\eta \in \mathbb{{R}^{N\times K}_{+}}} {\mathrm{min}} \qquad & \sum_{n\in N} \sum_{k\in K}   \eta_{nk} \cdot \parallel \zeta_n - \xi_k  \parallel^2 
\\ & subject \ to
\\ & \sum_{n\in N} \eta_{nk} = p_{k} \qquad \forall k \in K=\{1,...,K\}
\\ & \sum_{k\in K} \eta_{nk} = p^{\prime}_{n} \qquad \forall n \in N=\{1,...,N\}
\\ & \sum_{n\in N} \sum_{k\in K} \eta_{nk} = \sum_{n\in N} p^{\prime}_{n} =  \sum_{k\in K} p_{k} = 1
\end{split}
\end{align}
where, $\eta_{nk}$ is the amount of probability shipped from realisations $\zeta_n$ of an original discrete distribution of $N$ samples with probabilities $p^{\prime}_{n}$ to the final reduced discrete distribution of $K$ scenarios with values $\xi_k$ and probabilities $p_k$ \citep{Bertsimas2022}.

\citet{Li2014} formulated the probabilistic distance problem as MILP (introduced by name OSCAR) which aims to the minimisation not only of the probabilistic Kantorovich distance between the original and final scenarios but also the differences of the expected values on the output of the problem. The latter concept is referred to as Output Space System Response, and is imposed in objective function to enforce minimisation between the best, worst and expected performance. This concept was extended in a follow-up work to a sequential setting in order to efficiently reduce computational costs in problems with a large number of uncertain parameters \citep{Li2016-SR}. Moreover, \citet{Li2016-IA} based on the transportation distance minimisation problem presented an LP-based iterative scenario reduction algorithm showcasing its lower computational complexity while efficiency on solving case studies is maintained. Recently, \citet{Kammammettu2023} used Sikhorn distance for scenario reduction proposing an MINLP model and an iterative algorithm that is indicated to achieve computational cost reduction and competent stability over different runs and an increasing number of scenario sets' sizes.

As mentioned, while a majority of scenario generation approaches pay attention only on the distribution or the statistical properties of the uncertain set, OSCAR by \citet{Li2014} can also integrate the performance of the output of the problem in its formulation. In a similar context, the problem-based scenario generation and reduction constitute a topical issue of research \citep{Keutchayan2021,Henrion2022}. To this end, \citet{Bertsimas2022} presented a novel optimisation based approach for scenario reduction in which they introduced the term "problem-dependent divergence". Minimising this quantity in variants of well-known scenario reduction algorithms, enhanced efficiency is demonstrated. Beyond optimisation-based methods, \citet{Silvente2019} proposed a problem-based scenario reduction methodology using sensitivity analysis and evaluation on the associated problem under study.

\subsection{Stability evaluation of scenario generation methods}
\label{sec:evaluation}

The quality of SG methods is determinant of the effective solution of SP problems. Hence, measures to evaluate the quality of SG methods have been proposed in the literature. Critical is the work by \citet{Kaut2007} who proposed a methodology to evaluate the stability of SG methods. Moreover, the work by \citet{Bayraksan2006} presents measures to assess the quality of stochastic solutions in stochastic programs, which generally complement the quality assessment of SG approaches.

\subsubsection{Optimal cost of the stochastic program}

Towards the presentation of the stability measures, the notion of the optimal cost of a stochastic program is elucidated. Here, the true stochastic process is denoted as $\xi$ and then the optimal objective value of the stochastic program is ${z}^{\ast}$. Neglecting the symbols regarding the second-stage variables, Eqs. \eqref{eq:tssp} \& \eqref{eq:expcost} can be compactly written as \citep{Bayraksan2006}:

\begin{align}\label{eq:ev1}
z^{\ast} = \underset{{x \in X}} {\mathrm{min}}  \ f\left( x ; {\xi} \right)
\end{align}

Stochastic programs may render unsolvable when $\xi$ represents a true continuous stochastic program or a discrete distribution with a very large number of realisations. In such a case, the stochastic process $\xi$ can be approximated by a large reference tree, $R$. However, it is noted that $R$ must be generated by an unbiased sampling method \citep{Kaut2007}. Eventually $R$ contains a large number of independent and identically distributed (i.i.d.) realisations $K$, denoted as $\xi_1,...,\xi_k$. Thus, the stochastic program can be approximately reformulated as the scenario-based problem of Eq. \eqref{eq:ev2}, with optimal value $z_K^\ast$ \citep{Bayraksan2006}.

\begin{align}\label{eq:ev2}
z^\ast_K = \underset{{x \in X}} {\mathrm{min}}  \ f\left( x ; R \right) = \underset{{x \in X}} {\mathrm{min}}  \ \frac{1}{K} \sum_{k=1}^{K} f\left( x ; \xi_k \right)
\end{align}

\subsubsection{In-sample and out-of-sample stability}

Stability measures are introduced by \citet{Kaut2007} and exclusively concern the quality testing of SG approaches. In particular, good stability results indicate the capability of a SG method to generate different trees which lead to consistent objective values. Considering multiple scenario trees by a certain SG method, $T_c$, which lead to corresponding stochastic solutions ${\bar{x}}_c$ of the scenario-based problem, the in-sample and out-of-sample stability are achieved if Eqs. \eqref{eq:iss} -\eqref{eq:oss} are true, respectively \citep{Kaut2007}:

\begin{align}\label{eq:iss}
f\left({\bar{x}}_c;T_c\right)\approx f\left({\bar{x}}_{c^\prime};T_{c^\prime}\right)
\end{align}
\begin{align}\label{eq:oss}
f\left({\bar{x}}_c;\xi\right)\approx f\left({\bar{x}}_{c^\prime};\xi\right)
\end{align}

As mentioned, a large reference tree $R$ can be used to approximate the stochastic process $\xi$ in Eq. \eqref{eq:oss} towards the out-of-sample stability assessment. Of particular interest is the stability evaluation in the so-called optimisation-based SG methods which are considered to lead to a unique scenario tree for a desirable size. In this case stability assessment over scenario trees of “slightly different sizes” (i.e., varying number of generated scenarios) is proposed in \citep{King2012}.

\subsubsection{Bias}

Apart from stability assessment for SG methods, the quality assessment of the stochastic solutions obtained is also indicative of their actual performance on the stochastic programs at hand. The goal is to identify if a stochastic solution obtained using the scenario-based problems induce approximation error to the expected solution of the problem using the true stochastic process. The stochastic solution, $\bar{x}$, obtained by the scenario-based problem using a scenario tree $T$, may introduce a bias to the solution of the true stochastic program that is quantified by Eq. \eqref{eq:bias1}. However, following the notion of previous sections, the true stochastic process $\xi$ may be approximated by a large reference tree $R$, and the bias is estimated as in Eq. \eqref{eq:bias2} \citep{Kaut2007}:

\begin{align}\label{eq:bias1}
B\left(\bar{x}\right)=f\left(\bar{x};\xi\right) - \underset{{x \in X}} {\mathrm{min}}  \ f\left( x ; \xi \right) =f\left(\bar{x};\xi\right)-z^\ast
\end{align}
\begin{align}\label{eq:bias2}
B\left(\bar{x}\right)=f\left(\bar{x};R\right) - \underset{{x \in X}} {\mathrm{min}}  \ f\left( x ; R \right) =f\left(\bar{x};R\right)-z^\ast_K
\end{align}

Intuitively, the stochastic solution of the scenario-based problem is of perfect quality if it leads to the same objective value as the true stochastic program. It is noted that bias may be also referred to as “optimality gap” or “approximation error" in the literature \citep{Bayraksan2006,King2012}. For this study the term “bias" is preferred to avoid confusion with measures of optimality gap using replication procedures \citep{Bayraksan2006}, or generally the optimality gap of MILP models' solution.

\subsection{Contribution of this work}
\label{sec:contribution}

This work focuses on the Distribution and Moment Matching Problem and aims to investigate their stability while proposing an enhanced modelling approach. MMP by \citet{Kaut2021} and DMP by \citet{Bounitsis2022} differ from previous MMP models as these generate scenarios from an original set of historical data. Intuitively, distribution and moment matching statistical measures have been transformed to the objective function terms of a scenario reduction model instead of the traditional probabilistic distances. However, the presence of multiple terms in the objective function may impose numerical issues depending on the selection of user-defined weights and so their stability has to be evaluated. The novel proposed game theoretic modelling approach of DMP aims to introduce an alternative notion to handle such optimisation-based scenario reduction problems with multiple terms of the objective function in the future.   

The majority of the aforementioned works on scenario generation and reduction methods evaluate the stability of the methods over a number of generated sets in order to validate their performance. However, this can not be the case for optimisation-based scenario reduction models which can be considered deterministic, in the notion that they can provide an unique set when are optimally solved. The work by \citet{Bounitsis2022} has evaluated stability testing over scenario sets of “slightly different sizes” demonstrating that relatively low errors are obtained for DMP MILP over varying final number of scenarios (especially compared to the errors of DMP NLP models). These results can be interpreted as enhanced out-of-sample stability of DMP MILP.

However, in this context the work by \citet{Bounitsis2022} demonstrates that the selection of the user-defined weights for the errors of the objective function of the SR model can impose variations on the expected values and the bias. Nonetheless, previous works on the DMP \citep{Calfa2014} and MMP MILP version \citet{Kaut2021} specify a certain set of weights for the solution of the optimisation-based SR models. In particular, \citet{Calfa2014} define weights to capture the percentage errors on the terms of the objective function. On the other hand, \citet{Kaut2021} selects arbitrary weights based on the fact that stochastic programs are typically more sensitive to errors in means and variance than for the higher moments \citep{Seljom2021,Chopra1993}. Despite the wealth of the literature, in all of the aforementioned works stability over the selection of these user-defined parameters is not studied.

In this work, the ultimate goal is to evaluate the stability of the existing DMP optimisation-based models over the user-defined weights, as previous work have shown sensitivity of the models. So the stability of the DMP NLP and MILP models is to be validated. Furthermore, a game theoretic optimisation extension of the DMP MILP model by \citet{Bounitsis2022} is presented, in which the different terms of the objective function of the scenario reduction model are considered as the players of the game. The proposed model remains MILP by employing a separable programming approach \citep{Gjerdrum2001}. Finally, stability of the SR models is assessed over a number of errors' weights sets. As DMP NLP is shown to display under-specification issues, which has been resolved by the DMP MILP version, the correlation between the alleviation of under-specification issues and increased efficiency of the SG method is studied.

In the following section, the scenario generation methodology for the extended game theoretic DMP model is detailed. 

\section{Methodology}
\label{sec:methodology}

\subsection{Mathematical formulation of Distribution \& Moment Matching Problem}
\label{sec:prelim}

In Section \ref{sec:methodology} the game theoretic extension of the DMP MILP model by \citet{Bounitsis2022} is presented. Original DMP MILP aims at scenario reduction through an errors' minimisation problem. DMP is modelled as MILP by exploiting the existence of a large set of original scenarios and the employment of the Manhattan distance ($L^{1}$-norm) or the Chebyshev distance ($L^{\infty}$-norm) to express the objective function. As the game theoretic model uses the aforementioned work as basis the first four statistical properties are to be matched, namely: mean, variance, skewness and kurtosis of the known distribution of an uncertain parameter. Especially, skewness adds valuable information regarding the risk asymmetry and finally kurtosis describes the thickness of the tails and points out the significance of extreme scenarios for a distribution \citep{King2012}. 

In order to define the first four statistical moments of a known distribution, sets $i \in I$, as the index for the uncertain parameters, and $n \in N = \{ 1,...,N \} $, as the index of original scenarios, are defined. Moreover, variable $x_{in}$ denotes the value of the uncertain parameter $i$ in the node $n$ and variable $p_n$ its corresponding probability of occurrence. First four statistical moments are calculated as in Eqs. \eqref{eq:mean}-\eqref{eq:kurt}:

\begin{align}\label{eq:mean}
{Mean}_{i} = \mu_{i} = \sum_{n \in N} x_{in} \cdot p_{n} \qquad \forall i \in I
\end{align}
\begin{align}\label{eq:var}
{Variance}_{i} = \sigma_{i}^{2} = \sum_{n \in N} ( x_{in} - \mu_i )^{2} \cdot p_{n} \qquad \forall i \in I
\end{align}
\begin{align}\label{eq:skew}
{Skewness}_{i} = \frac{\sum_{n \in N} ( x_{in} - \mu_i )^{3} \cdot p_{n}} {\sigma_{i}^{3}}  \qquad \forall i \in I
\end{align}
\begin{align}\label{eq:kurt}
{Kurtosis}_{i} = \frac{\sum_{n \in N} ( x_{in} - \mu_i )^{4} \cdot p_{n}} {\sigma_{i}^{4}}  \qquad \forall i \in I
\end{align}

It is noted that skewness and kurtosis are by definition normalised properties, but for the sake of MMP definition these are denormalised. By this option, nonlinearities due to scaling are avoided but deviations of mean of variance in the final set can impose even more intense deviations for the errors of skewness and kurtosis \citep{Calfa2014,Kaut2021}. 

\citet{Calfa2014} were the first to introduce the distribution matching extension to MMP. Especially, generated sets are enforced to additionally match the marginal Empirical Cumulative Distribution Function (ECDF) of the uncertain parameters. \citet{Calfa2014} incorporated ECDF in the model through approximation by the Generalised Logistic Function (GLF) in a pre-process step. In contrast, \citet{Bounitsis2022} directly match ECDF values via the DMP MILP model. Theoretically, Cumulative Distribution Function (CDF) expresses the probability of a random parameter, $z_i$, to take a value lower or equal to some value $t$. Empirical CDF (ECDF) is a non-parametric estimator of the CDF which is defined by Eq. \eqref{eq:ecdf}:

\begin{align}\label{eq:ecdf}
ECDF\left(t\right) = \frac{1}{n} \cdot \sum_{i=1}^{N} \textbf{1}\{z_i \leq t\}
\end{align}
where, $n$ is the sample size and $\mathbf{1}\{z_i\le t\}$ the indicator function.

\subsection{Nash game theoretic approach}
\label{sec:nash}

Nash game theoretic approach is the cornerstone of the reformulation of the original DMP MILP model in this study. For the sake of brevity the interested reader is referred to \citet{Marousi2023} for a detailed presentation of game-theoretic approaches, while a brief outline is presented as follows. Nash bargaining approach proposes that taking into consideration the initial position (status quo) of all players of a game, then a fair solution among all players can be achieved \citep{Nash1950}. Let us denote a profit maximisation problem for each player of a game with $t \in T$ players. Given the status quo solutions of the players $t \in T$ prior to joining the game (or generally a lower profit requirement point), $\pi_t^{SQ}$, then using a Nash game theoretic approach a fair solution point, $\pi_t$, can be achieved. The Nash bargaining solution approach is able to fairly contribute the payoff to the players obeying the axioms of: (i) Pareto optimality, (ii) symmetry, (iii) linear invariance, and (iv) independence of irrelevant alternatives. For a typical profit maximisation problem, the fair solution aims to maximise the Nash product, $\Psi$, which is given in Eq. \eqref{eq:nash1} \citep{Harsanyi1977}:

\begin{align}\label{eq:nash1}
\Psi = \prod_{t \in T} \left( \pi_t - \pi_t^{SQ} \right)^{\alpha_t}  
\end{align}
where $\alpha_t$ is the negotiation power of each player $t$. Hence, a player $t$ enters the game only if $\pi_t \geq \pi_t^{SQ}$, i.e., if the profit entering the game is greater than its status quo point. 

The objective function of Eq. \eqref{eq:nash1} is nonlinear and nonconvex leading to computationally challenging problems. However, recently various reformulation strategies have been proposed in the PSE literature towards the efficient approximation of the Nash product \citep{Marousi2023}. This work adopts the logarithmic transformation coupled with SOS2 variables for piecewise approximation of Nash product in order to alleviate nonlinearities and retain the MILP formulation of the DMP model by \citet{Bounitsis2022}. This reformulation strategy achieves the linearisation of the Nash product based on a separable programming approach and has been initially proposed by \citet{Gjerdrum2001}. Firstly, a convexification step is conducted by employing a logarithmic transformation on the Nash product. Then the concave Nash product is linearised via a piecewise linear function of $g \in G$ prespecified grid points. So, the linear approximation of the Nash product, $\overline{\Psi}$, is computed as in Eq. \eqref{eq:nash2}:

\begin{align}\label{eq:nash2}
\overline{\Psi} = \sum_{t \in T} \sum_{g \in G} \alpha_t \cdot ln\left( \overline{\pi}_{tg} - \pi_t^{SQ} \right)  \cdot \lambda_{tg}
\end{align}
where, $\lambda_{tg}$ are the SOS2 variables for the computation of the profit of player $t$ at grid point $g$, which is denoted as $\overline{\pi}_{tg}$. Eventually the profit of each player $t$, $\pi_t$, is estimated by the combination of Eqs. \eqref{eq:nash3} \& \eqref{eq:nash4}:

\begin{align}\label{eq:nash3}
\sum_{g \in G} \lambda_{tg} = 1 \qquad \forall t \in T 
\end{align}

\begin{align}\label{eq:nash4}
\pi_t = \sum_{g \in G} \overline{\pi}_{tg} \cdot \lambda_{tg} \qquad \forall t \in T 
\end{align}

\subsection{DMP via game theoretic approach}
\label{sec:model}

\begin{figure}[H] 
\centering
\includegraphics[width=16.4cm]{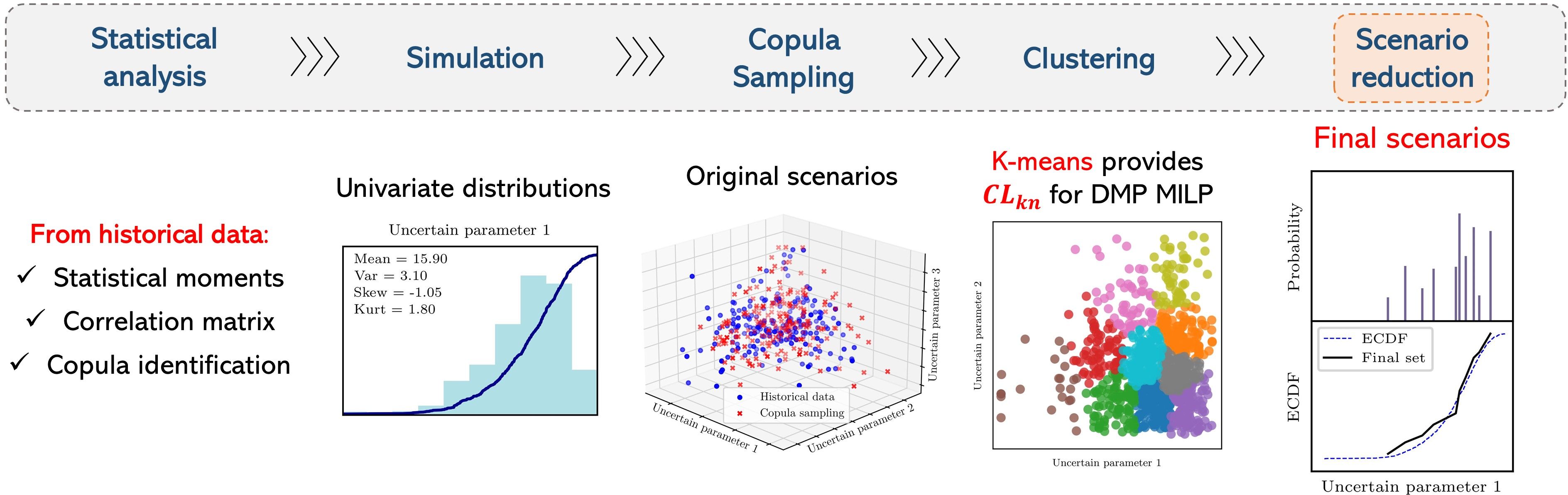}
\caption{Proposed extended DMP MILP optimisation model can be incorporated in the data-driven scenario generation framework presented by \citet{Bounitsis2022}. }
\label{fig:dmpmeth}
\end{figure}

In this section the extended game theoretic version of DMP MILP model for scenario reduction is presented. This model can replace the original DMP MILP on the data-driven scenario generation methodology by \citet{Bounitsis2022}. For the sake of completeness the framework is briefly outlined. Using historical data as input, the statistical properties regarding the moments ($D_{im}$), correlations ($C_{ii^{\prime}}$) and copula are initially estimated/constructed. Then using the results of the statistical analysis, univariate distributions $U_{in}$ are simulated and then copula-based original scenarios $X_{in}$ can be generated. In the next step, original scenarios $X_{in}$ are clustered into as many clusters as the desirable size of the reduced set $K$. Hence, we introduce a dynamic set ${CL}_{kn}$, which maps each scenario $n \in N$ to exactly one prespecified cluster $k \in K$ using the labels of clustering. From the original scenarios only one is selected at each prespecified cluster $k \in K$, through binary variables $y_{kn}$, and probabilities of occurrence $p_{kn}$ are defined. A visualisation of the steps is provided in Fig. \ref{fig:dmpmeth}.  

As for the DMP MILP, for the presented game theoretic models only $L^1$-norm or $L^{\infty}$-norm are used to quantify the errors. Especially, the game theoretic model remains MILP by employing the separable programming approach for Nash product approximation outlined in Section \ref{sec:nash}. The model is presented for a game of 3 players which are introduced by index $t \in T = \{ SM,ECDF,COV \}$. These players are the errors on moments, covariance and ECDF matching which also correspond to the terms of the objective function of the original model (in particular, for ECDF the maximum for each parameter is captured, $e_{i}$). The mathematical model using $L^1$-norm for the quantification of errors is presented as follows: 

\begin{align}\label{eq:objnash}
 \underset{{y_{kn},p_{kn}}} {\mathrm{max}} \qquad \sum_{t \in T} \sum_{g \in G} ln\left( \pi_t^{MAX} - \overline{\pi}_{tg} \right)  \cdot \lambda_{tg}
\end{align}
\begin{align*}
subject\ to
\end{align*}
\begin{align}\label{eq:mnash1}
\sum_{g \in G} \lambda_{tg} = 1 \qquad \forall t \in T 
\end{align}
\begin{align}\label{eq:mnash2}
\pi_t = \sum_{g \in G} \overline{\pi}_{tg} \cdot \lambda_{tg} \qquad \forall t \in T 
\end{align}
\begin{align}\label{eq:mnash3}
\pi_{'SM'} = \sum_{i \in I} \sum_{m \in M} W^{SM}_{m} \cdot ( d^{+}_{im} + d^{-}_{im} ) 
\end{align}
\begin{align}\label{eq:mnash4}
\pi_{'COV'} = \sum_{i \in I} \sum_{i^{\prime} \in I,i^{\prime}>i} W^{COV}_{ii^{\prime}} \cdot ( c^{+}_{ii^{\prime}} + c^{-}_{ii^{\prime}} ) 
\end{align}
\begin{align}\label{eq:mnash5}
\pi_{'ECDF'} = \sum_{i \in I} W^{ECDF}_{i} \cdot e_{i}
\end{align}
\begin{align}\label{eq:y1}
\sum_{n:(k,n) \in CL_{kn}} y_{kn} = 1 \qquad \forall k \in K
\end{align}
\begin{align}\label{eq:y2}
\sum_{k:(k,n) \in CL_{kn}} y_{kn} \leq 1 \qquad \forall n \in N
\end{align}
\begin{align}\label{eq:p1}
P^{min} \cdot y_{kn} \leq p_{kn} \leq P^{max} \cdot y_{kn}  \qquad \forall (k,n) \in CL_{kn}
\end{align}
\begin{align}\label{eq:p2}
\sum_{(k,n) \in CL_{kn}} p_{kn} = 1
\end{align}
\begin{align}\label{eq:m1}
\sum_{(k,n) \in CL_{kn}} \left( X_{in} \cdot p_{kn} \right) + d^{+}_{im} - d^{-}_{im} = D_{im} \qquad \forall i \in I, m=1
\end{align}
\begin{align}\label{eq:m2}
\sum_{(k,n) \in CL_{kn}} \left[ {\left( X_{in} - D_{i,1} \right)}^{m} \cdot p_{kn} \right] + d^{+}_{im} - d^{-}_{im} = D_{im} \qquad \forall i \in I, m>1
\end{align}
\begin{align}\label{eq:c1}
\sum_{(k,n) \in CL_{kn}} \left[ \left( X_{in} - D_{i,1} \right)  \cdot \left( X_{i^{\prime}n} - D_{i^{\prime},1} \right) \cdot p_{kn} \right] + c^{+}_{ii^{\prime}} - c^{-}_{ii^{\prime}} = C_{ii^{\prime}} \qquad \forall i,i^{\prime},i<i^{\prime}
\end{align}
\begin{align}\label{eq:e1}
\sum_{k \in K} y_{kn} \cdot ECDF_{in} - \sum_{(k,n^{\prime}) \in CL_{kn^{\prime}} \bigwedge X_{in^{\prime}} \leq X_{in} }  p_{kn^{\prime}} = \phi_{in} \qquad \forall i \in I, n \in N
\end{align}
\begin{align}\label{eq:e2}
e_{i} \geq \phi_{in} - \left(  1 - \sum_{k:(k,n) \in CL_{kn}} y_{kn} \right)   \qquad \forall i \in I, n \in N
\end{align}
\begin{align}\label{eq:e3}
e_{i} \geq - \phi_{in} - \left(  1 - \sum_{k:(k,n) \in CL_{kn}} y_{kn} \right)   \qquad \forall i \in I, n \in N
\end{align}
\begin{align}\label{eq:b1}
-1 \leq \phi_{in} \leq 1  \qquad \forall i \in I, n \in N
\end{align}
\begin{align}\label{eq:b2}
c^{+}_{ii^{\prime}}, c^{-}_{ii^{\prime}} \geq 0  \qquad \forall i \in I, i^{\prime} \in I
\end{align}
\begin{align}\label{eq:b3}
d^{+}_{im}, d^{-}_{im} \geq 0  \qquad \forall i \in I, m \in M
\end{align}
\begin{align}\label{eq:b4}
e_{i} \geq 0  \qquad \forall i \in I
\end{align}
\begin{align}\label{eq:b5}
0 \leq p_{kn} \leq 1  \qquad \forall (k,n) \in CL_{kn}
\end{align}

Eq. \eqref{eq:objnash} constitutes the objective function of the model and is combined with \eqref{eq:mnash1} - \eqref{eq:mnash5} to introduce the game theoretic optimisation approach for the 3 players $t \in T$. Eq. \eqref{eq:y1} impose the selection of exactly one original scenario $n$ from each cluster $k$. Eq. \eqref{eq:y2} complements Eq. \eqref{eq:y1} imposing that an original scenario $n$ must be selected for at most one cluster. Eq. \eqref{eq:y2} is redundant but it is added to the model as enhances performance and reduces the number of nodes needed. Eq. \eqref{eq:p1} defines the allowable limits for the selected scenarios (selection through binary variable $y_{kn}$), and Eq. \eqref{eq:p2} enforces the cumulative probability of all final scenarios to be equal to 1. Eqs. \eqref{eq:m1} - \eqref{eq:m2} estimate the errors regarding the first four moments between the original and the reduced distributions. Analogously, Eq. \eqref{eq:c1} estimates the errors regarding the correlation matrix in cases that multiple uncertain parameters $i$ are taken into consideration. Eq. \eqref{eq:e1} estimates the deviations, $\phi_{i,n}$, regarding the ECDF curve between the original and the reduced distribution for every final scenario $n$ and each parameter $i$. Then, Eqs. \eqref{eq:e2} - \eqref{eq:e3} capture the maximum error, $e_{i}$, regarding the ECDF matching for each parameter $i$ over only the selected scenarios $n$. Finally, Eqs. \eqref{eq:b1} - \eqref{eq:b5} set the variables' bounds. The combination of the variables for the various types of errors and the problem's objective function, Eq. \eqref{eq:objnash}, leads to their efficient minimisation.

For the definition of the model for $L^{\infty}$-norm the only change pertains to the quantification of the errors. Thus, \eqref{eq:mnash3} - \eqref{eq:mnash5} are simply replaced by equations:

\begin{align}\label{eq:mnash6}
\pi_{'SM'} \geq W^{SM}_{m} \cdot ( d^{+}_{im} + d^{-}_{im} ) \qquad \forall i \in I, m \in M
\end{align}
\begin{align}\label{eq:mnash7}
\pi_{'COV'} \geq  W^{COV}_{ii^{\prime}} \cdot ( c^{+}_{ii^{\prime}} + c^{-}_{ii^{\prime}} ) \qquad \forall i \in I, i^{\prime} \in I
\end{align}
\begin{align}\label{eq:mnash8}
\pi_{'ECDF'} \geq W^{ECDF}_{i} \cdot e_{i} \qquad \forall i \in I
\end{align}

In particular, the objective function of Eq. \eqref{eq:objnash} aims at the maximisation of the Nash product, which for this case is expressed using as status-quo solutions, the parameters $\pi_t^{MAX}$, which represent the higher allowable statistical errors for each term $t \in T$. So, status-quo $\pi_t^{MAX}$ can be captured by solving the original DMP MILP model from \citet{Bounitsis2022} omitting the errors regarding term $t$ from the objective function. The error value regarding the omitted term $t$ after the solution of the model is considered as the maximum allowable error (status quo). Hence, an additional step to SG framework is necessary when the game theoretic DMP is implemented.

For the calculation of errors' values crucial is the role of the user-defined values of errors' weights. A relatively larger value for a weight indicate that the corresponding attribute is forced to be matched more accurately. Weights of the errors regarding the statistical moments and the covariance depend on user-defined parameters $\overline{W}^{SM}_{im}$ and $\overline{W}^{COV}_{ii^{\prime}}$, respectively. So, weights are defined as shown in Eqs. \eqref{eq:w1} - \eqref{eq:w2}.

\begin{align}\label{eq:w1}
{W}^{SM}_{im} = \frac{\overline{W}^{SM}_{im}}{\mid D_{im} \mid} \qquad \forall i \in I, m \in M
\end{align}
\begin{align}\label{eq:w2}
{W}^{COV}_{ii^{\prime}} = \frac{\overline{W}^{COV}_{ii^{\prime}}}{\mid C_{ii^{\prime}} \mid} \qquad \forall i \in I, i^{\prime} \in I
\end{align}

Furthermore, $W^{ECDF}_{i}$ is the weight regarding the error of the ECDF matching and its high values enforce accurate matching of the ECDF curve. Denoting as $\overline{W}^{ECDF}_{i}$ the user-defined parameter regarding the ECDF error weight and considering that ECDF values range from 0 to 1, it is set that $W^{ECDF}_{i}$ is equal to $\overline{W}^{ECDF}_{i}$ ($W^{ECDF}_{i} = \overline{W}^{ECDF}_{i}$ ).

Last but not least, a crucial point for the solution of MMP and DMP MILP models regards the exact matching of the mean through Eq. \eqref{eq:m1}. In particular, \eqref{eq:m1} estimates the mean value from the selected nodes of the original set but the computed value may deviate from the mean of the original distribution. A deviation of the mean of the reduced set may be problematic for the MMP and DMP MILP problems because rest errors on moments and correlation matrix use the mean of the original distribution through Eqs. \eqref{eq:m2} \& \eqref{eq:c1}, respectively. In order to avoid interdependent approximate estimations regarding the latter statistical properties, Eq. \eqref{eq:malter} is proposed to replace Eq. \eqref{eq:m1} and so mean is enforced to not deviate from its actual value.

\begin{align}\label{eq:malter}
\sum_{(k,n) \in CL_{kn}} X_{in} \cdot p_{kn}  = D_{im} \qquad \forall i \in I, m=1
\end{align}

Nonetheless, this condition may lead to infeasibility during the solution of the DMP MILP when the original set of scenarios is relatively small or more uncertain parameters are considered \citep{Xu2012,Kaut2021}.

Ultimately, the proposed framework uses historical data as input and its steps are summarised in Algorithm \ref{alg:framework}. 

\RestyleAlgo{ruled}
\begin{algorithm}[H] 
\caption{Data-driven SG framework incorporatin game theoretic DMP MILP.} \label{alg:framework}
\textbf{Input:}\ Historical data regarding uncertain set of $\mid I \mid$ parameters; \\
\textbf{Step 1:}\ Statistical analysis $\rightarrow$  Moments $D_{im}$, Covariance $C_{ii^{\prime}}$, Copula $\Psi$; \\
\textbf{Step 2:}\ Simulation of a univariate Pearson distributions based on $D_{im}$ \ $\rightarrow$  $U_{in}$ ; \\
\textbf{Step 3:}\ Copula-based sampling $\rightarrow$ $\mid N \mid$ original scenarios $X_{in}$ : \\
\uIf{$  \mid I \mid = 1  $}{
   Identical to the Pearson distribution of \textbf{Step 2} \;
   }
\uElseIf{$  \mid I \mid \geq 2  $}{
   Sampling pair- or Vine- copula $\Psi$ \ $\rightarrow$ Pseudo observations $O_{in}$ ; \\
   Linear interpolation of $O_{in}$ to $U_{in}$ \ $\rightarrow$ original set $X_{in}$ ;
   }
\textbf{Step 4:}\ K-means clustering of $X_{in}$ to $\mid K \mid$ clusters $\rightarrow$ Labels on subsets $CL_{kn}$ \\
\textbf{Step 5:}\ DMP MILP omitting a player $t$ at a time $\rightarrow$ Maximum allowable errors $\pi_t^{MAX}$  \\
\textbf{Step 6:}\ Game theoretic DMP MILP $\rightarrow$ Selections $y_{kn}$ \& probabilities $p_{kn}$ of $\mid K \mid$ final scenarios; \\
\textbf{Output:}\ Set of $\mid K \mid$ generated scenarios with values $y_{kn} \cdot X_{in} $ \& probabilities $p_{kn}$; \\
\end{algorithm}

\section{Case study: Capacity planning under uncertainty}
\label{sec:casestudy}

\subsection{Summary of the problem}
\label{sec:casestudyintro}

\begin{figure}[H]  
\centering
\includegraphics[width=10.4cm]{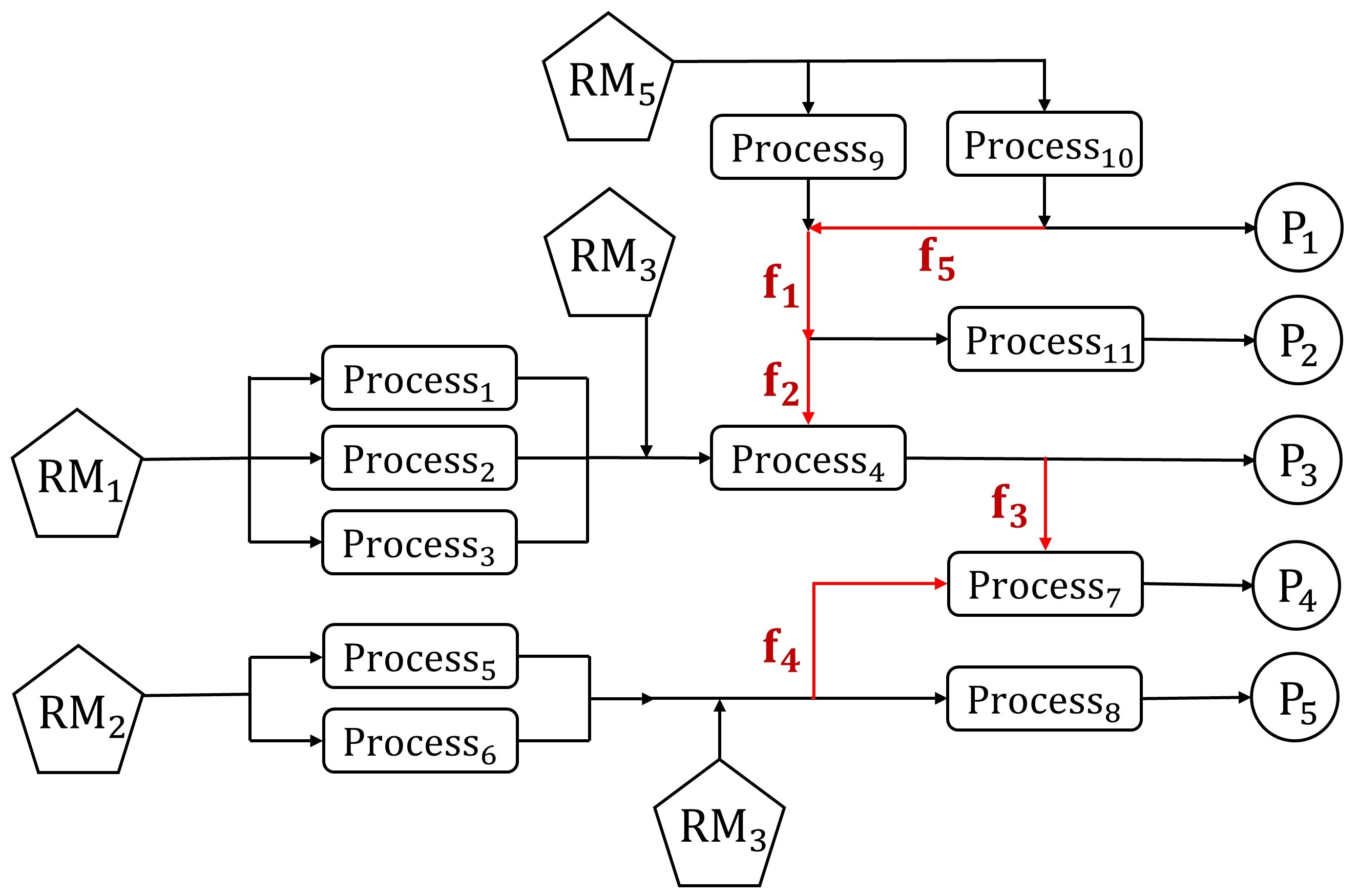}
\caption{Process network for the investigated case study.}
\label{fig:csn}
\end{figure}

A capacity planning MILP problem under uncertainty is studied. It is originally examined by \citet{Acevedo1998}, and several case studies of the stochastic problem are presented in \citep{Li2014}. It considers 5 products ($l \in L$) that are produced using 5 raw materials ($h \in H$) and 11 candidate processes ($j \in J$). The flowsheet of the process is visualised in Figure \ref{fig:csn}.

The optimal decisions regard the selection of the candidate processes to be installed through binary variables, $y_j$. The deterministic version of the MILP problem is formulated according to flowsheet in Fig. \ref{fig:csn} and is modelled as follows:

\begin{align} \label{eq:cs2obj}
\underset{{y_j}} {\mathrm{max}} \qquad \sum_{l=1}^{5} \beta_{l} \cdot P_{l} - \sum_{h=1}^{5} \alpha_{h} \cdot {RM}_h - \sum_{j=1}^{11} {OC}_j \cdot {IS}_j - \sum_{j=1}^{11} \left( {DC}_j \cdot Q_j + {FC}_j \cdot y_j \right)
\end{align}
\begin{align*}
subject \ to
\end{align*}
\begin{align} \label{eq:cs2:1}
{OS}_j = {PC}_j \cdot IS_j \qquad \forall j=\{1,..,11\}
\end{align}
\begin{align} \label{eq:cs2:2}
P_l = D_l \qquad \forall l=\{1,..,5\}
\end{align}
\begin{align} \label{eq:cs2:3}
{RM}_h \leq {maxRM}_h \qquad \forall h=\{1,..,5\}
\end{align}
\begin{align} \label{eq:cs2:4}
{IS}_j \leq {MI}_j \cdot Q_j \qquad \forall j=\{1,..,11\}
\end{align}
\begin{align} \label{eq:cs2:5}
Q_j \leq {maxQ}_j \cdot y_j   \qquad \forall j=\{1,..,11\}
\end{align}
\begin{align} \label{eq:cs2:6}
{RM}_1={IS}_1+{IS}_2+{IS}_3
\end{align}
\begin{align} \label{eq:cs2:7}
{RM}_2={IS}_5+{IS}_6
\end{align}
\begin{align} \label{eq:cs2:8}
{RM}_3+{OS}_1+{OS}_2+{OS}_3+f_2={IS}_4
\end{align}
\begin{align} \label{eq:cs2:9}
f_3+f_4={IS}_7
\end{align}
\begin{align} \label{eq:cs2:10}
{RM}_5={IS}_9+{IS}_{10}
\end{align}
\begin{align} \label{eq:cs2:11}
{OS}_1+f_5=f_1
\end{align}
\begin{align} \label{eq:cs2:12}
{OS}_{10}=f_5+P_1
\end{align}
\begin{align} \label{eq:cs2:13}
f_1=f_2+{IS}_{11}
\end{align}
\begin{align} \label{eq:cs2:14}
{OS}_4=f_3+P_3
\end{align}
\begin{align} \label{eq:cs2:15}
P_2={OS}_{11}
\end{align}
\begin{align} \label{eq:cs2:16}
P_4={OS}_7
\end{align}
\begin{align} \label{eq:cs2:17}
P_5={OS}_8
\end{align}

The description of the model is given as follows:

\begin{itemize}
    \setlength\itemsep{-1pt}
    \item eq. \eqref{eq:cs2obj} is the objective function which aims at profit maximisation,
    \item constraint \eqref{eq:cs2:1} defines the yield relations,
    \item constraint \eqref{eq:cs2:2} enforces demands for the products,
    \item constraint \eqref{eq:cs2:3} enforces raw materials availability limit,
    \item constraints \eqref{eq:cs2:4} – \eqref{eq:cs2:5} indicate materials' flow in the network,
    \item constraints \eqref{eq:cs2:6} – \eqref{eq:cs2:17} define mass balances.
\end{itemize}

Uncertainty is assumed to be realised in the yield constant of processes $j$ (${PC}_j$), in the recourse matrix, and in the demand of products $l$ (${D}_l$), on the right hand side (RHS). The problem under uncertainty can be modelled as TSSP and scenarios to describe the uncertainty for the uncertain parameters can be introduced. The selection of processes, $y_j$, and their capacities, $Q_j$, are set as here-and-now decisions, while all rest variables are wait-and-see decisions. In \ref{sec:appendix} the nomenclature and the values of the parameters regarding the case study is provided.

The performance of 3 scenario generation models is assessed:
\begin{itemize}
    \setlength\itemsep{-1pt}
    \item DMP NLP models by \citet{Calfa2014} - denoted as "NLP",
    \item DMP MILP models by \citet{Bounitsis2022} - denoted as "MILP",
    \item proposed Nash game theoretic reformulation of this work - denoted as "NASH",
\end{itemize}

All studied models are defined for both the $L^1$ and the $L^{\infty}$ norms. Moreover, MILP-based models uses copula-based original scenarios as input which are generated according to Algorithm \ref{alg:framework}. These sets of original scenarios are also used as reference set for the estimation of bias (for all three models). The expected solution of the TSSP for this benchmark set (tree $R$ of the representation in \ref{sec:evaluation}) is also denoted as full-space (FS) solution.

Especially, the evaluation of bias and stability for the investigated methods is conducted for a number of set of varying values regarding the errors' weights in the objective function. These set are constructed for the values of $\overline{W}^{SM}_{im}$, $\overline{W}^{COV}_{ii^{\prime}}$, $\overline{W}^{ECDF}_{i}$ as the investigated models are formulated considering 3 terms on the objective function, i.e., the players of the game. Typical values 1, 10, 50 and, 100 are considered and 37 unique sets of values are constructed. Regarding Nash reformulation through the separable programming approach, 50 grid points ($g$) are employed for the piece-wise linear approximations. Finally, three instances of increasing complexity are examined for the considered case study. These are summarised in Table \ref{table:cs-instances}.

\begin{table}[H]
\renewcommand{\arraystretch}{1.5}
\fontsize{10pt}{10pt} \selectfont
\begin{center}
\caption{ Investigated instances of the case study. }
\label{table:cs-instances}
\begin{tabular}{l c c c c } 

\hline
{ } &
\begin{tabular}{@{}c@{}} \textbf{\# uncertain} \\ \textbf{parameters}\end{tabular} & 
\begin{tabular}{@{}c@{}} \textbf{Uncertain} \\ 
\textbf{parameters}\end{tabular} & 
\begin{tabular}{@{}c@{}} \textbf{Original} \\ 
\textbf{scenarios}\end{tabular} & 
\begin{tabular}{@{}c@{}} \textbf{Final} \\ 
\textbf{Scenarios}\end{tabular}
\\	
\hline
\textbf{Instance 1} & 2 &	${PC}_7,\ {PC}_8$ &	1,000 &	10 \\ 
\textbf{Instance 2} & 4 &	${PC}_4,\ {PC}_7,\ {PC}_8,\ {PC}_{11}$ &	2,000 &	20 \\ 
\textbf{Instance 3} & 8 &	\begin{tabular}{@{}c@{}} ${PC}_4,\ {PC}_7,\ {PC}_8,\ {PC}_{11}$ \\ 
$D_2,\ D_3,\ D_4,\ D_5$ \end{tabular} &	2,000 &	40 \\ 
\hline

\end{tabular}
\end{center}
\end{table}

The executions are performed in a Dell workstation with Intel\textregistered \ Core\textsuperscript{TM} i9-10900K CPU @ 3.70 GHz and 32.00 GB RAM. The DMP NLP model \citep{Calfa2014} is solved using the BARON solver \citep{Sahinidis1996}, while MILP models using the GUROBI 9.5 solver within GAMS 38.2 modelling system \citep{gurobi,GAMS}. Although calculation of lower bound constitutes a significant drawback of DMP MILP, the Nash-based reformulations can achieve near optimal solutions in short execution times. Towards a fair comparison, time limit of 900s and an optimality gap tolerance of 1\% are set for the MILP models.

\begin{figure}[H] %
    \centering
    \subfloat[\centering Uncertain distribution of parameter $PC_7$ in Instance 1.] 
    {\label{fig:udist1} \includegraphics[width=7.5cm]{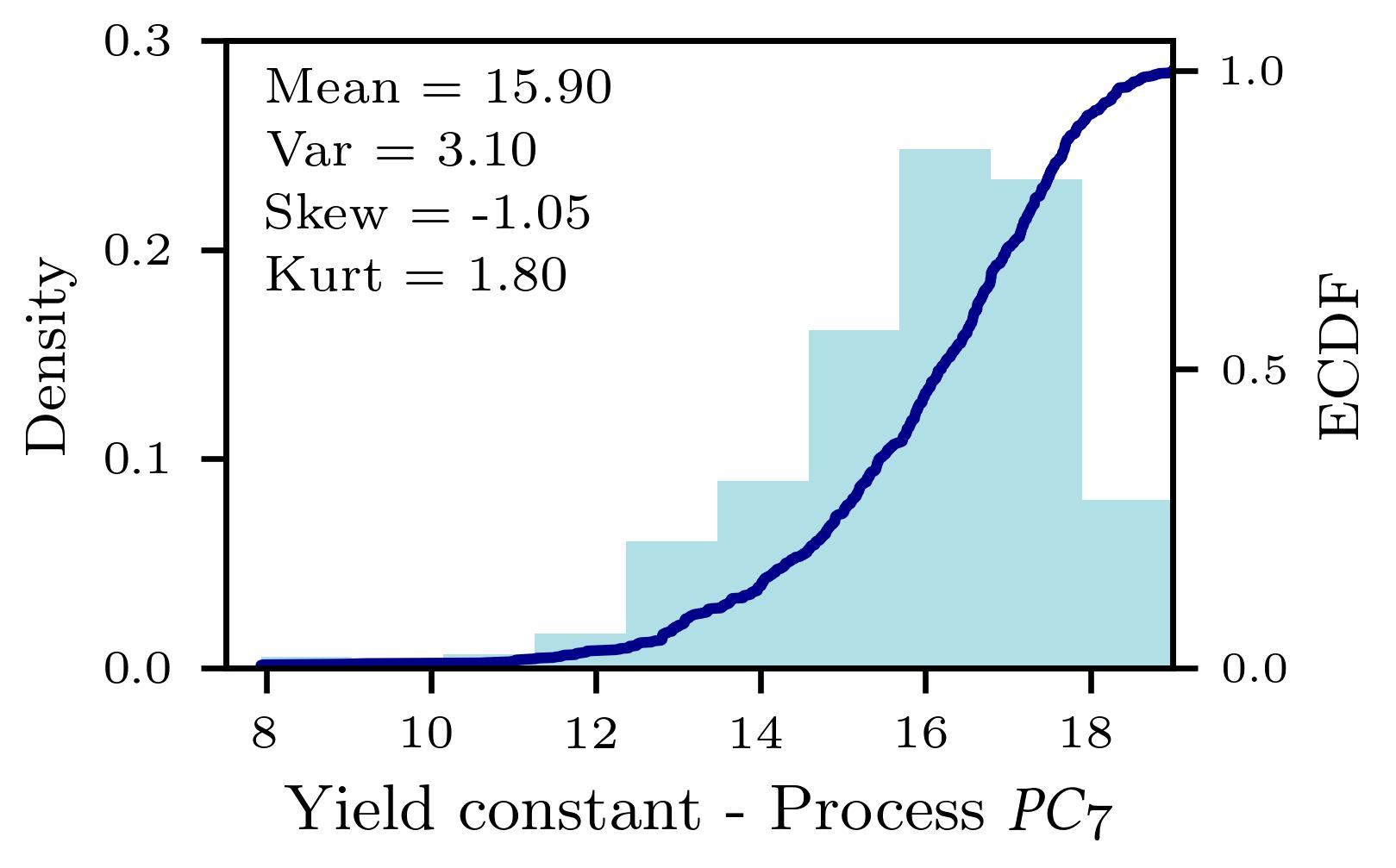} }%
    \quad
    \subfloat[\centering Uncertain distribution of parameter $PC_8$ in Instance 1.] 
    {\label{fig:udist2} \includegraphics[width=7.5cm]{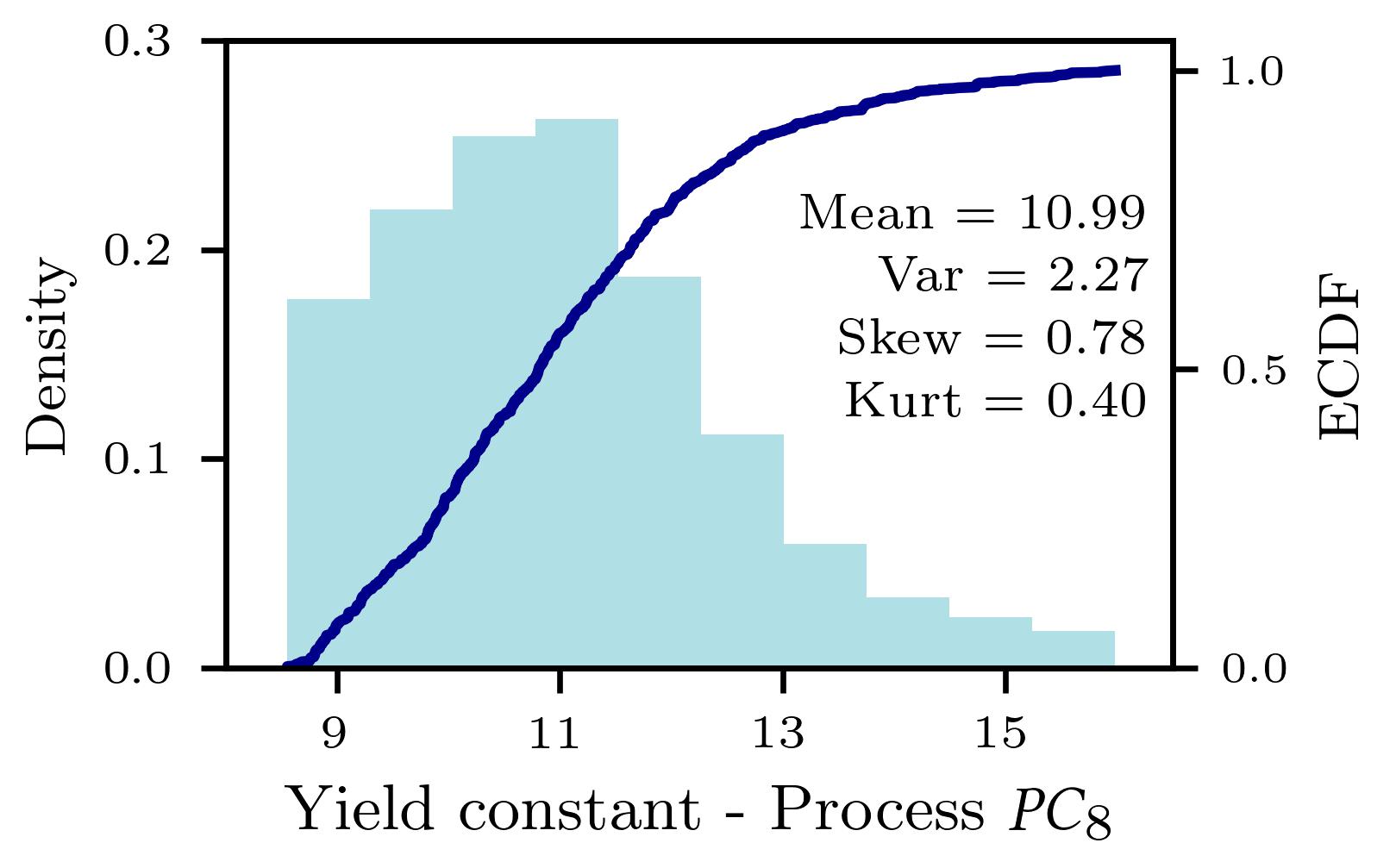} }%
    \\
    \vspace{12pt}
    \subfloat[\centering Historical data and copula-based sampled original scenarios for the bivariate case of Instance 1.] 
    {\label{fig:copulacomp} \includegraphics[width=7.5cm]{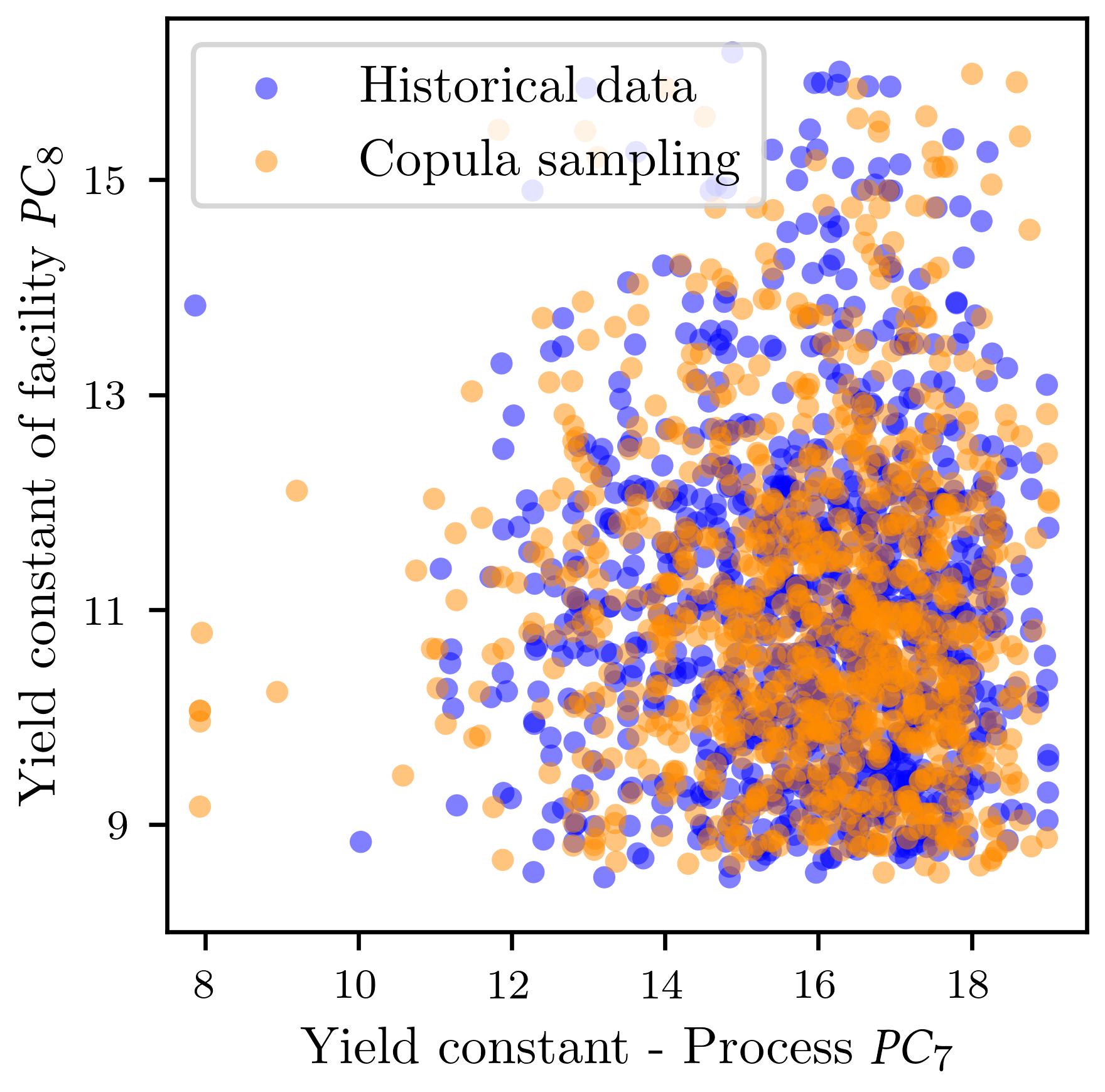} }%
    \subfloat[\centering Clustering of 1,000 original scenarios to 10 clusters for the bivariate case of Instance 1.] 
    {\label{fig:clusters} \includegraphics[width=7.5cm]{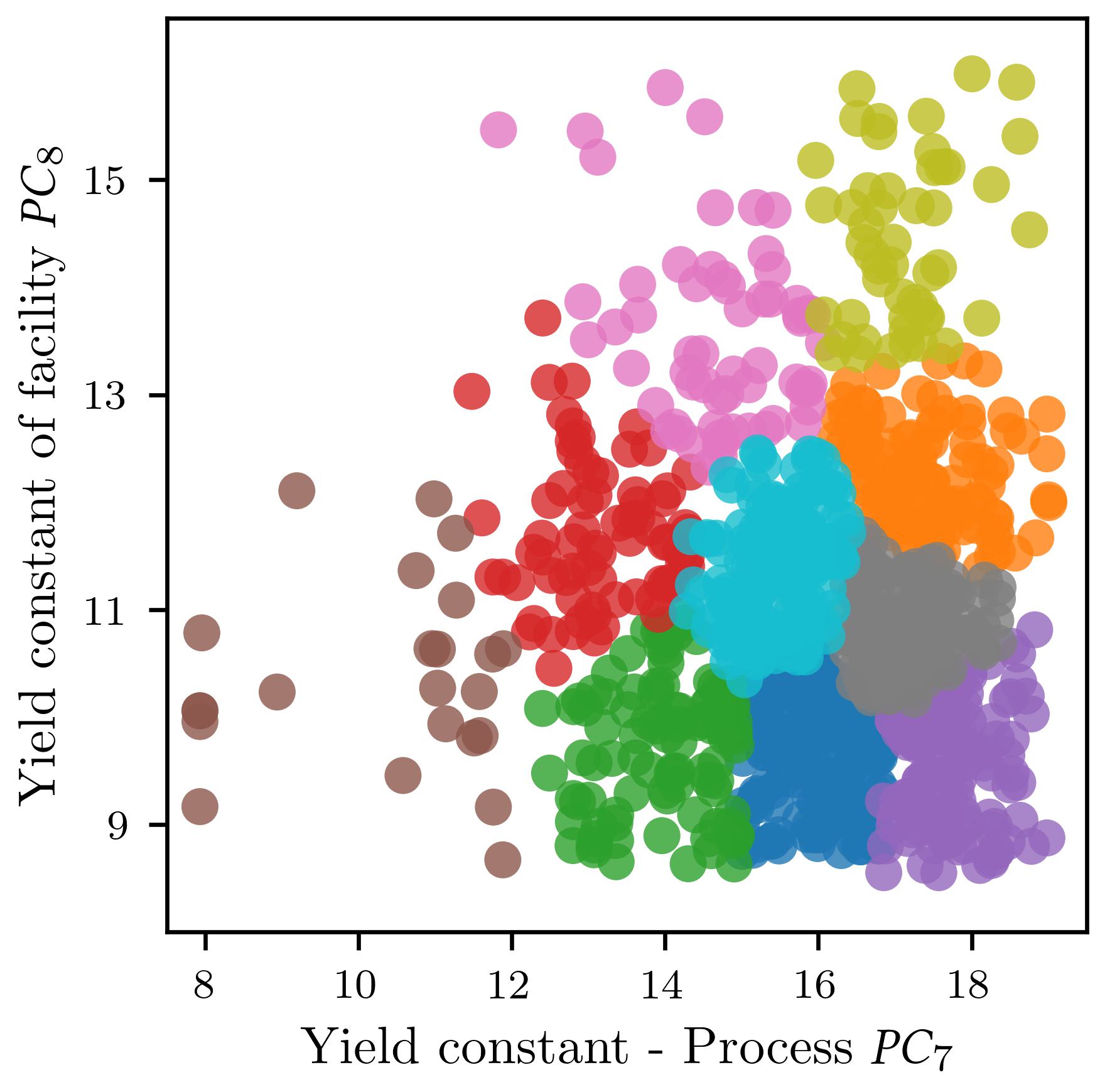} }%
    \qquad
    \caption{Initial steps (1 to 4) of the proposed scenario generation framework of Algorithm \ref{alg:framework}.}%
    \label{fig:original}%
\end{figure}

Before the discussion of results, some visualisations on the preliminary methodologies to generate original scenarios (input of DMP MILP models) utilising Algorithm \ref{alg:framework} are provided in Fig. \ref{fig:original}. These visualisations examine the bivariate case of Instance 1. The good match of the copula sampled original scenarios to historical data is of particular interest. More results regarding the scalability of the methods used in Algorithm \ref{alg:framework} can be found in \citet{Bounitsis2022}. 

\subsection{Quality of stochastic solutions}
\label{sec:quality}

In the work by \citet{Bounitsis2022} the scalability of the framework's techniques and the stability of the DMP MILP model over an increasing number of generated scenarios has been demonstrated. However, the latter work mainly focuses on quality assessment of statistical matching of the generated sets and on tests concerning the quality of the stochastic solutions. In this work the quality of stochastic solutions over the various sets of errors' weights is initially investigated. Hence, in Table \ref{table:bias} are presented bias results regarding all the considered scenario generation models. As 37 runs are executed for each model and each instance the minimum, maximum and mean values of bias are reported.

\begin{table}[H]
\renewcommand{\arraystretch}{1.5}
\fontsize{10pt}{10pt} \selectfont
\begin{center}
\caption{ Minimum, mean and maximum values of bias regarding all optimisation-based DMP models and instances over 37 sets of user-defined weights. }
\label{table:bias}
\begin{tabular}{l c c c c c c }

\hline
\% $\mid$ \textbf{Bias} $\mid$ & \textbf{NLP - $L^1$} & \textbf{NLP - $L^{\infty}$} & \textbf{MILP - $L^1$} & \textbf{MILP - $L^{\infty}$} & \textbf{NASH - $L^1$} & \textbf{NASH - $L^{\infty}$} 
\\	
\hline
{} & \multicolumn{6}{c}{\textit{Instance 1 – 2 uncertain parameters}}\\
\textbf{Min} &  0.112 & 0.109	 &	0.002 & 0.002	 & \textbf{0.001}	 & \textbf{0.000}	 \\ 
\textbf{Mean} & 1.327  & 1.383 & 0.047 & 0.034 & \textbf{0.027}	 & \textbf{0.028}	  \\
\textbf{Max} & 2.038  & 2.008	 & 0.171	 & 0.083	 & \textbf{0.073}	 & \textbf{0.076}	  \\
{} & \multicolumn{6}{c}{\textit{Instance 2 – 4 uncertain parameter}}\\
\textbf{Min} &  0.305 & 0.367	 & 0.011	 & 0.009	 & \textbf{0.006}	 & \textbf{0.002}	 \\ 
\textbf{Mean} & 2.364  & 2.834	 & 0.048	 & 0.049	 & \textbf{0.034}	 & \textbf{0.029}	  \\
\textbf{Max} & 3.349  & 3.430	 & 0.127	 & 0.208	 & \textbf{0.085}	 & \textbf{0.069}	  \\
{} & \multicolumn{6}{c}{\textit{Instance 3 – 8 uncertain parameter}}\\
\textbf{Min} &  0.739 &	2.412 &	\textbf{0.009} & \textbf{0.004}	 & 0.012	 & \textbf{0.004}	 \\ 
\textbf{Mean} & 2.761 &	3.049 &	0.058 &	\textbf{0.042} &	\textbf{0.042} &	\textbf{0.042}  \\
\textbf{Max} & 3.440  & 3.612	 & 0.210	 & 0.162	 & \textbf{0.087}	 & \textbf{0.082}	  \\
\hline
\end{tabular}
\end{center}
\end{table}

Results of Table \ref{table:bias} indicate that NLP models lead to sets whose stochastic solutions can impose big bias over the true stochastic problem (represented by a large reference set). The lowest achieved bias by NLP models are higher than the mean values of the MILP counterparts. Moreover, bias imposed by NLP SG models seems to increase significantly with the increase on the number of uncertain parameters. On the other hand, MILP models achieve much lower bias and the NASH reformulation even enhance the quality of the stochastic solutions. The achieved lower bias values are close to zero and especially the  $L^{\infty}$-norm-based version report slightly lower bias than the $L^1$-norm-based counterparts for each MILP model. Finally, NASH models significantly improve the values regarding the maximum bias compared to the MILP models. This result can indicate the mitigation of outlier values regarding the bias with the extension of MILP to a game theoretic approach. 

\begin{figure}[H] %
    \centering
    \subfloat[\centering Normalised bias over the sets of weights.] 
    {\label{fig:heat-bias} \includegraphics[width=12.5cm]{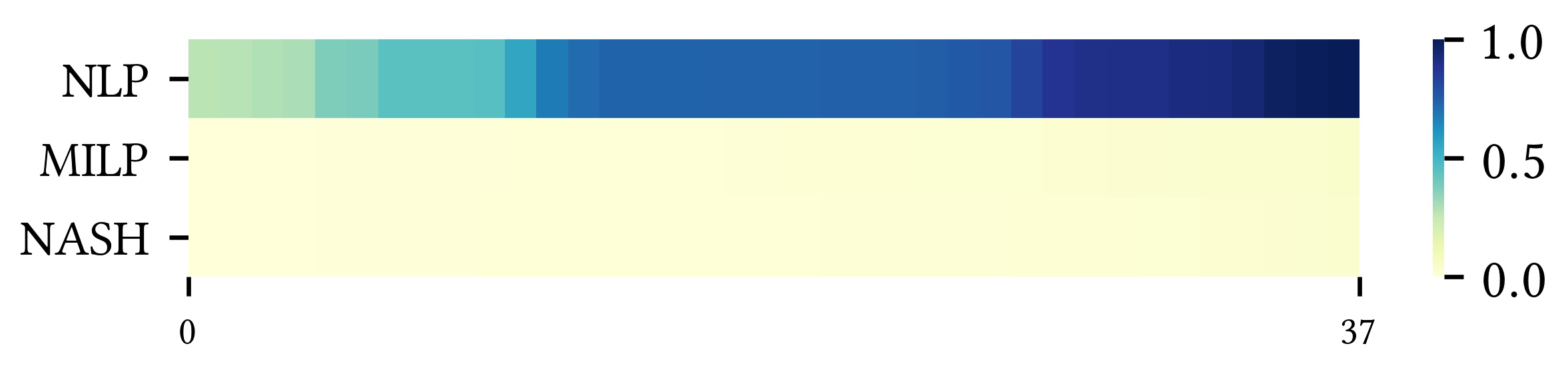} }%
    \\
    \vspace{12pt}
    \subfloat[\centering Normalised moment matching absolute maximum errors over the sets of weights.] 
    {\label{fig:heat-mm} \includegraphics[width=12.5cm]{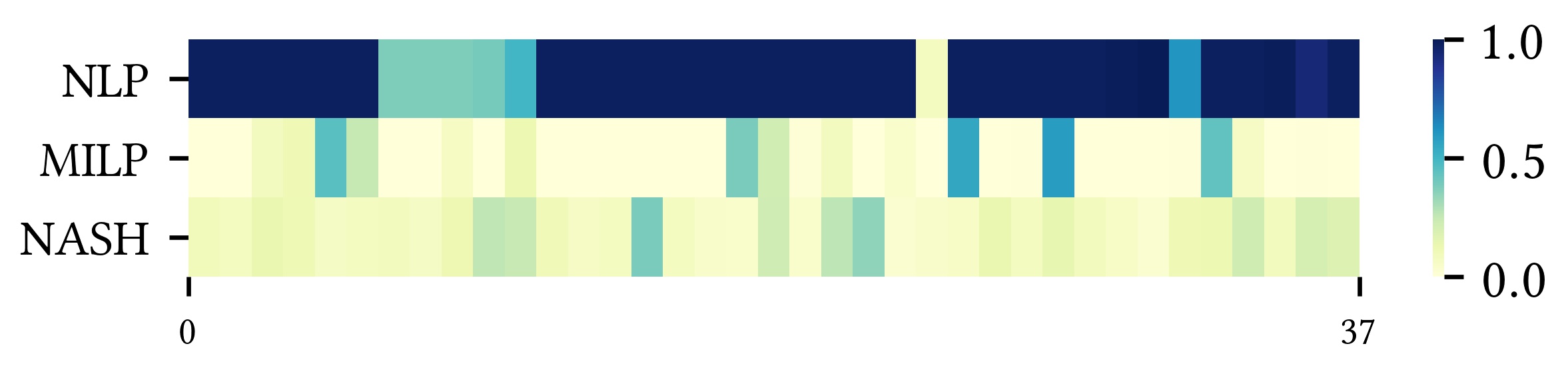} }%
    \\
    \vspace{12pt}
    \subfloat[\centering Normalised ECDF matching absolute maximum errors over the sets of weights.] 
    {\label{fig:heat-ecdf} \includegraphics[width=12.5cm]{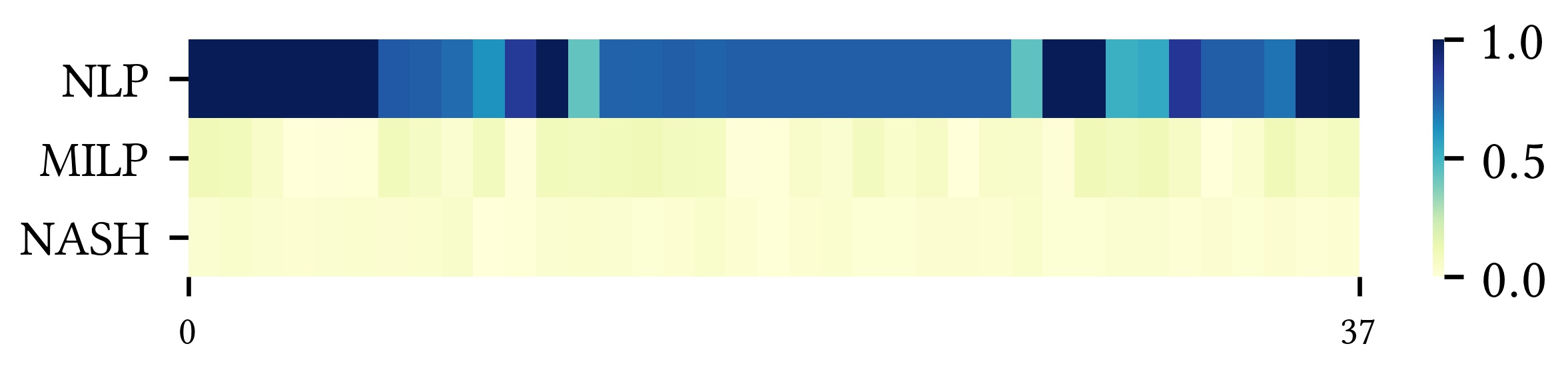} }%
    \caption{Bias, moments and ECDF matching normalised errors obtained using the $L^{\infty}$-norm-based versions of DMP models over the errors' weights sets in Instance 1. The data are normalised for each of the attribute and sorted from minimum to maximum values of bias for each SR method (NLP, MILP and NASH).}%
    \label{fig:heatmap}%
\end{figure}

A question that arises from the bias results is how much the errors on the statistical matching of the SG models can affect the quality of the stochastic solutions using the generated sets. Thus, in Fig. \ref{fig:heatmap}   correlation matrices in order to identify the correlation between the statistical errors and the bias for all the $L^{\infty}$-norm-based versions of the SG models of Instance 1 are presented. It is noted that in each plot the values are normalised (based on minimum and maximum values of the relevant measure) and sorted from the minimum to maximum values of bias of each SG model.

Fig. \ref{fig:heatmap} demonstrates that the high bias values regarding the NLP model are mostly owed to the order of magnitude of the statistical errors. Conversely, the low bias values of the MILP models merely come due to the relatively low values on the statistical matching errors and are less affected by the exact matching in the statistical sense. In other words, instances with higher moment matching or ECDF matching errors do not utterly lead to higher values of bias. Overall, MILP and NASH methods can constantly lead to scenario sets of good quality independently of cases with relatively higher errors on the statistical matching. Finally, MILP seems to perform slightly better regarding the moment matching will NASH outperforms on ECDF matching. 

A possible explanation of the latter results and the high bias of NLP models may come if the under-specification issues are taken into consideration. In particular, NLP models often display numerical issues during their executions which may lead to under-specified scenario sets. These issues are more intense for cases with more uncertain parameters. In this case study NLP models consistently lead to a lower number of generated scenarios than the prespecified ones. For Instance 1, 3-6 scenarios are generated from NLP models instead of the 10 prespecified scenarios. However, this small number of generated scenarios can not consistently capture the statistical properties of the original distributions. So, NLP models may be very helpful on the fast and effective construction of scenario trees as presented by \citet{Calfa2014}, but generally their under-specification issues are resolved by the MILP formulations, which consequently enhance the quality on scenario reduction on the investigated case study. 

\subsection{Statistical matching of SG models}
\label{sec:statistics}

\begin{figure}[H] %
    \centering
    \subfloat[\centering Moments and ECDF matching absolute maximum errors obtained using the $L^{\infty}$-norm-based versions of the models.] 
    {\label{fig:errors1} \includegraphics[width=7.5cm]{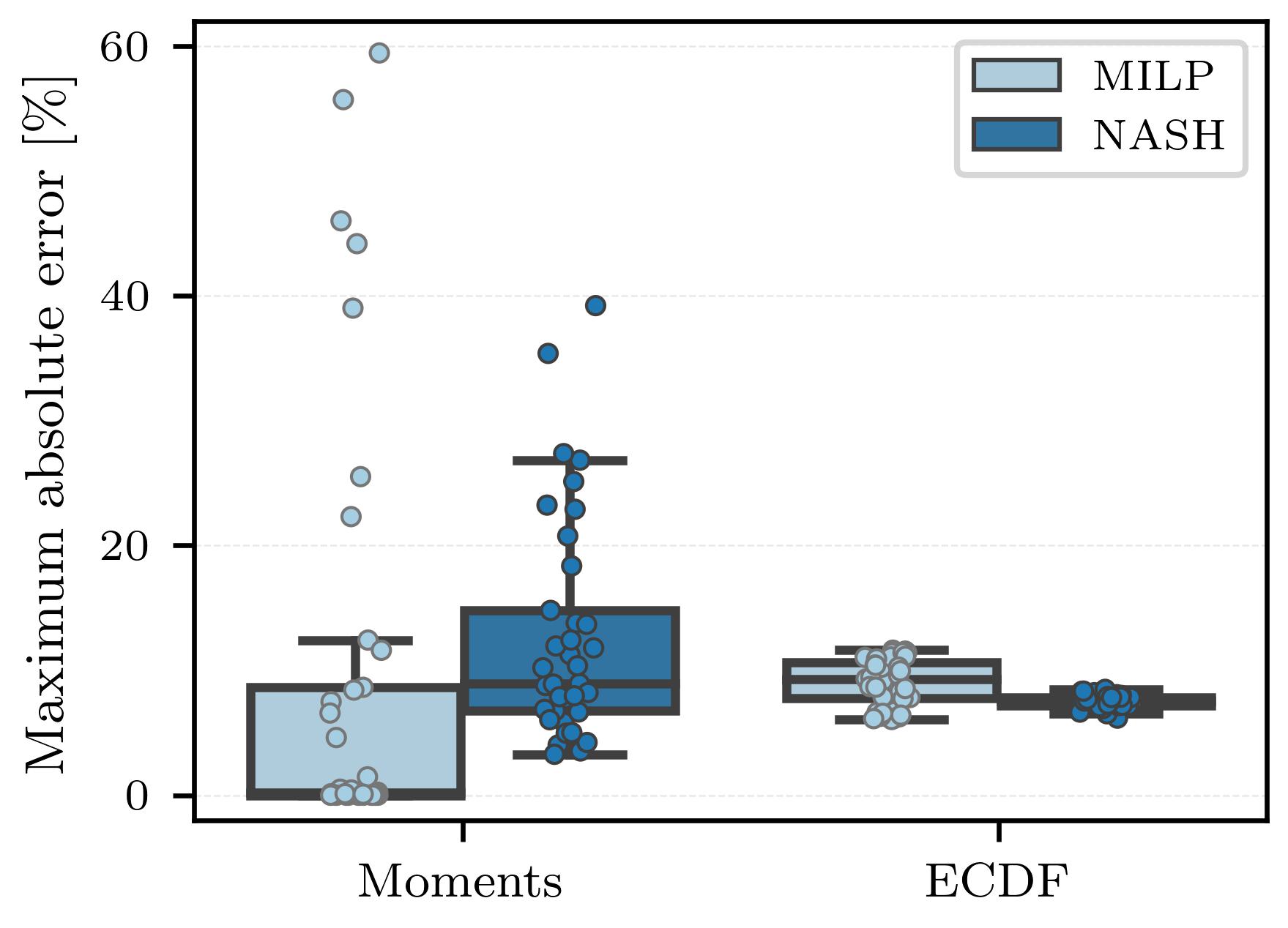} }%
    \\
    \vspace{12pt}
    \subfloat[\centering Moments matching absolute errors using the $L^{\infty}$-norm-based versions of the models.] 
    {\label{fig:errors2} \includegraphics[width=13.4cm]{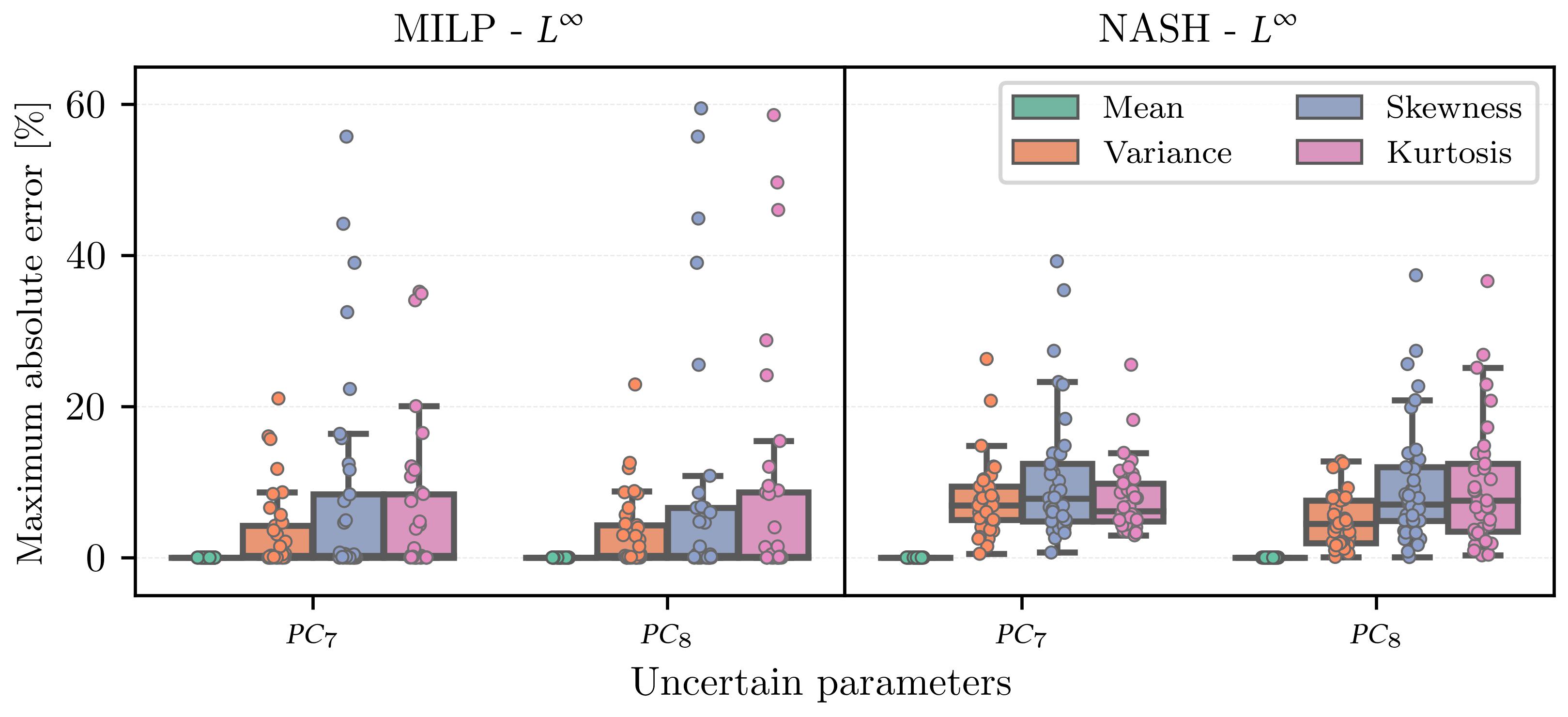} }%
    \caption{Statistical matching results regarding DMP MILP and NASH $L^{\infty}$-norm-based models in Instance 1 over the errors' weights sets.}%
    \label{fig:errors}%
\end{figure}

As demonstrated by Fig. \ref{fig:heatmap} NLP DMP models lead to errors of statistical properties which are multiple times bigger than the errors obtained by corresponding MILP models. Focusing on MILP and NASH DMP models for Instance 1, in Fig. \ref{fig:errors} results regarding the errors of the models on the statistical matching are reported. Especially, these results are reported for the $L^{\infty}$-norm-based versions of the models as the maximum absolute errors could be more comprehensible for such a discussion.

In Fig. \ref{fig:errors1} are initially reported the errors regarding the terms of the objective functions (or players in NASH model) and then in Fig. \ref{fig:errors2} are reported the errors on the first for statistical moments over the sets of errors' weights. First of all, it is noted that both NASH and MILP models achieve zero errors regarding the error on the covariance matrix for Instance 1 and so this is not envisaged. The results of Fig. \ref{fig:errors1} highlight the smaller variance of NASH generated scenario sets on the matching errors regarding the moments and ECDF terms of the objective function. Regarding the errors on the moments, although NASH model achieves generally higher errors compared to the original MILP, the variance of the errors is smaller and extreme outliers are avoided. Similarly, the extended NASH model leads to importantly lower variance on the maximum ECDF error, while the mean ECDF error over all sets is also lower than the original MILP method.

Regarding the specific moments' errors, detailed results are presented in Fig. \ref{fig:errors2}. As mentioned for the total moments matching error, NASH model lead to reduced variance, while the values of the specific errors is generally increased compared to MILP model. Nonetheless, the instances of the most extreme outliers in the errors of the sets using the MILP model are avoided by the sets generated from the proposed NASH version.

\begin{figure}[H] %
    \centering
    \subfloat[\centering Generated scenarios regarding $PC_{7}$.] 
    {\label{fig:scenarios-7} \includegraphics[width=15.5cm]{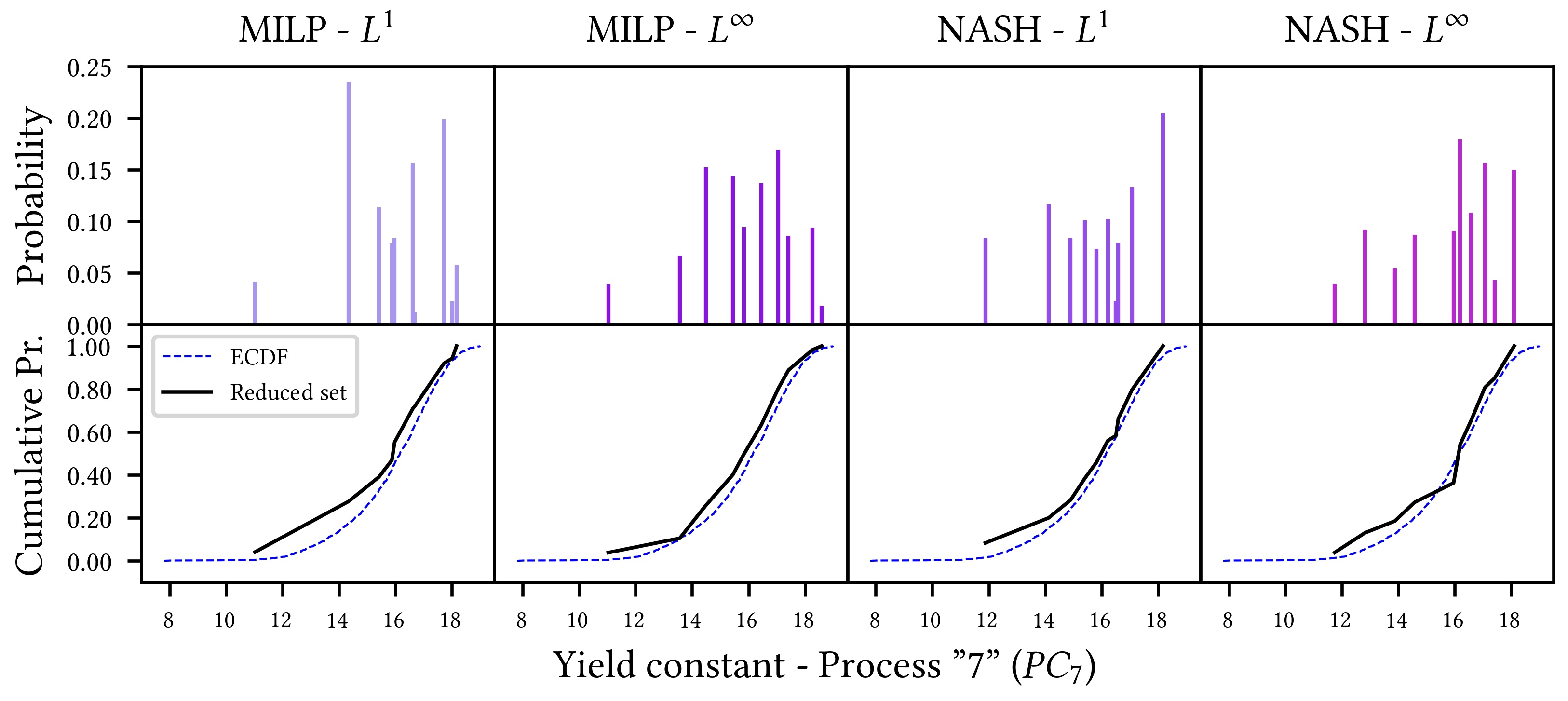} }%
    \\
    \vspace{12pt}
    \subfloat[\centering Generated scenarios regarding $PC_{8}$.] 
    {\label{fig:scenarios-8} \includegraphics[width=15.5cm]{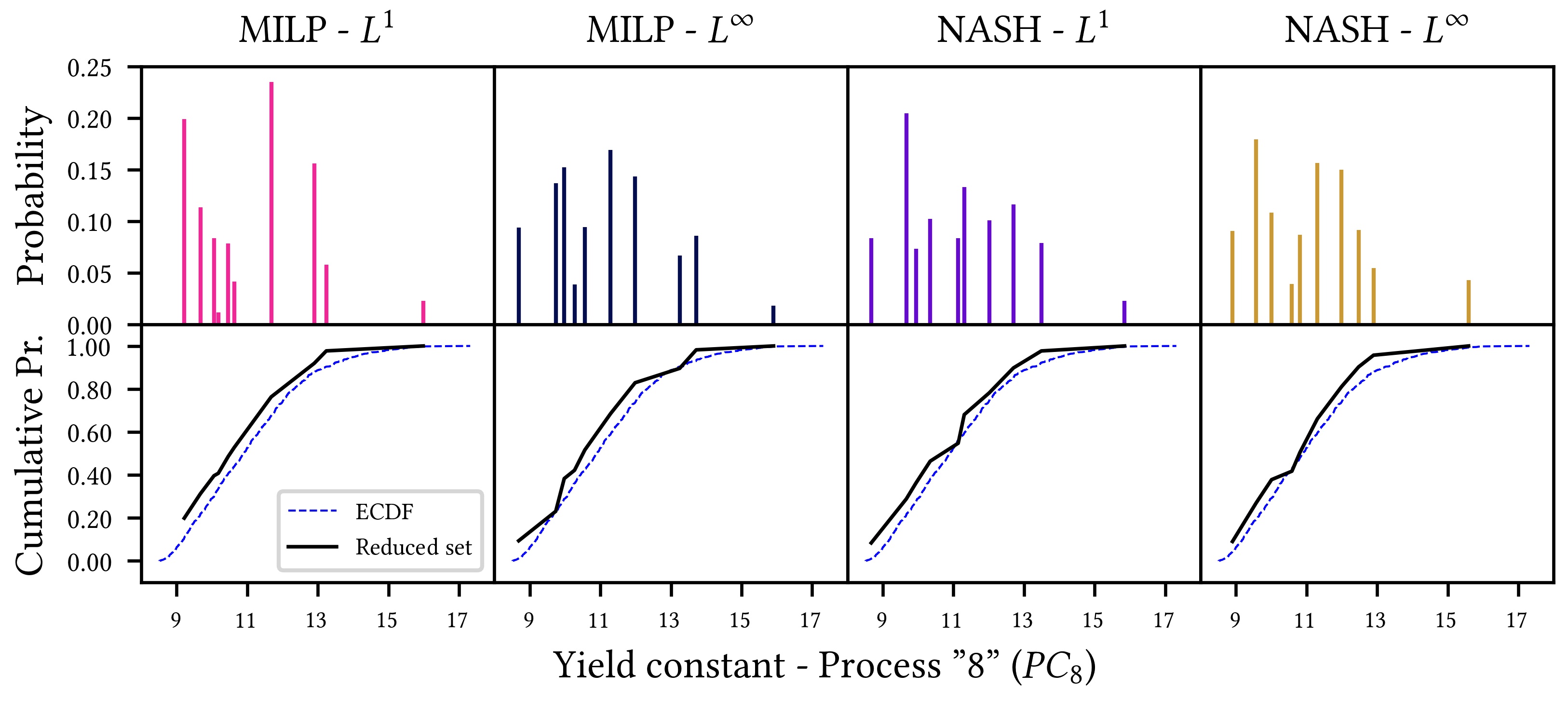} }%
    \caption{Values of scenarios and ECDF matching regarding $PC_{7},\ PC_{8}$ in Instance 1 for ($\overline{W}^{SM}_{im}$, $\overline{W}^{COV}_{ii^{\prime}}$, $\overline{W}^{ECDF}_{i}$) = (1,1,1).}%
    \label{fig:scenarios}%
\end{figure}

Finally, an indicative visualisation of the generated sets and their ECDF matching for Instance 1 for a certain set of weights' errors is visualised in Fig. \ref{fig:scenarios}. These visualisations indicate the good matching in the statistical sense of both MILP and NASH models. Even though only 10 scenarios are generated to capture the uncertainty on both uncertain parameters $PC_{7}$ \& $PC_{8}$, ECDF of the original data for each of the parameters is sufficiently approximated. Similarly, even the distribution of the values and their probabilities when visualised are comparable to the original distributions of Fig. \ref{fig:original}. Overall, the reduced sets approximate well in the statistical sense the original distribution and this is reflected also to the quality of the stochastic solutions.

\subsection{Stability of scenario generation methods}
\label{sec:stability}

As presented in Section \ref{sec:quality} the use of MILP and NASH DMP models can result to scenario sets that lead to stochastic solutions of good quality as assessed by the bias measure. Nevertheless, in Section \ref{sec:statistics} it is demonstrated that there is an accurate and consistent matching of the statistical properties of the DMP models over the 37 sets of the errors' weights. Finally, in this section is evaluated the stability of the MILP SG models in order to investigate whether the game theoretic extension of the DMP model can enhance its performance. In particular, in-sample and out-of-sample stability are evaluated according to their definitions in Section \ref{sec:evaluation} and using the sets of original scenarios as the reference tree (denoted by $R$) to represent the true stochastic process. The objective value of the TSSP using the reference tree is also noted as full-space (FS) solution. Ultimately, the effectiveness of the DMP models to consistently provide scenario trees that lead to expected values with lower divergence and good accuracy of the TSSP problem's solution at hand is evaluated.

\begin{figure}[H] %
    \centering
    \subfloat[\centering Instance 1: In-sample stability.] 
    {\label{fig:2p-iss} \includegraphics[width=8cm]{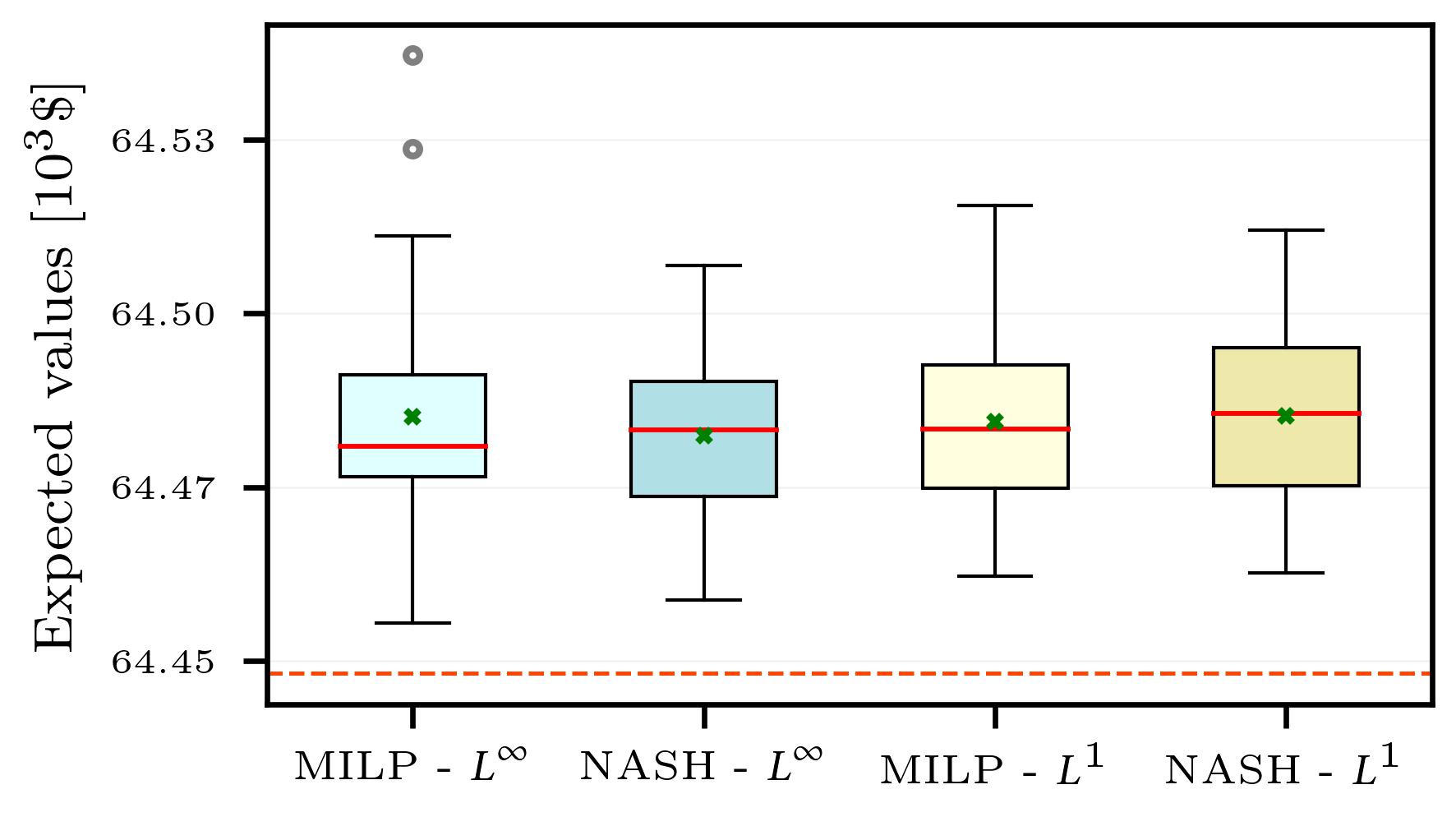} }%
    \subfloat[\centering Instance 1: Out-of-sample stability.] 
    {\label{fig:2p-oss} \includegraphics[width=8cm]{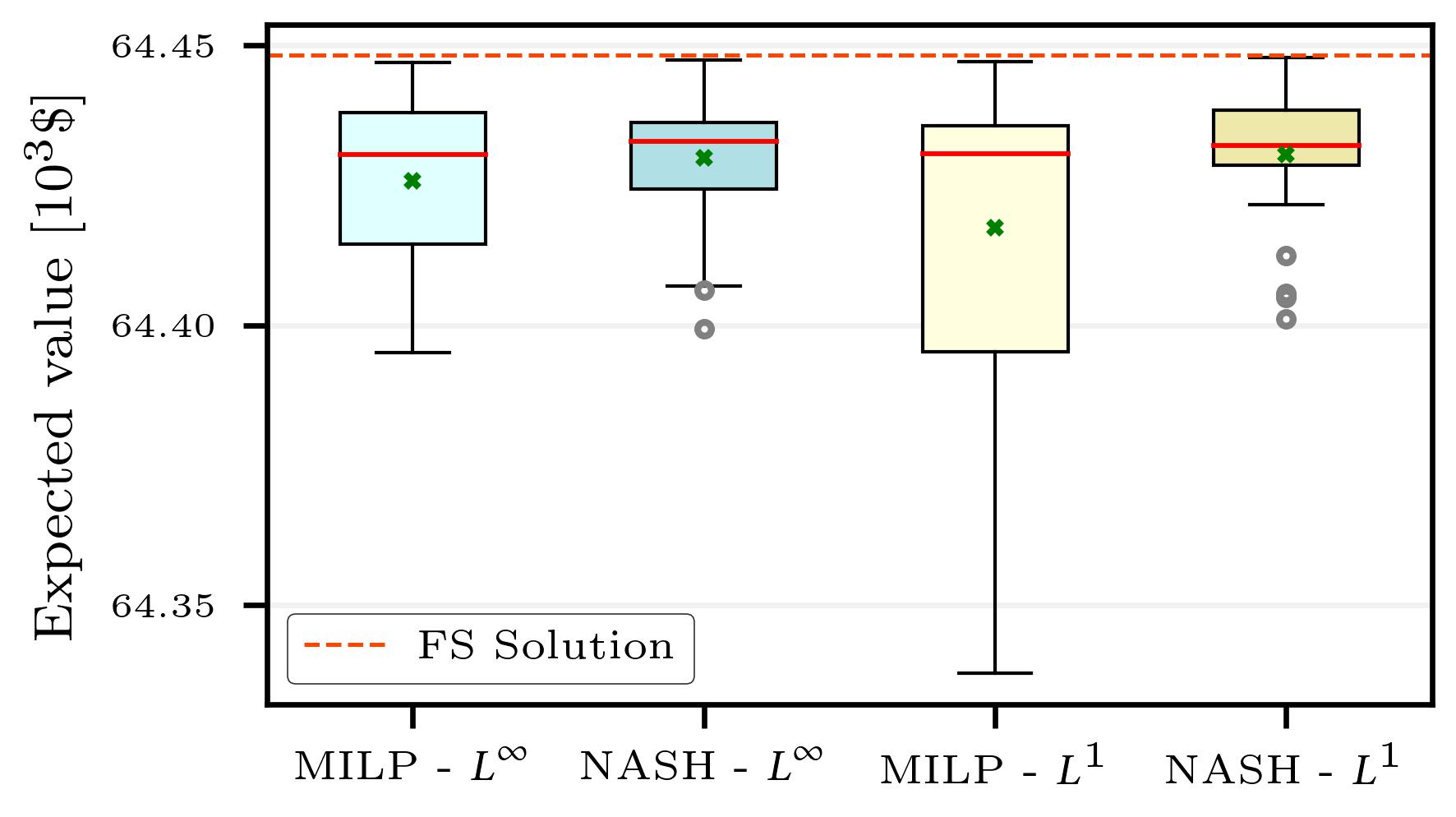} }%
    \\
    \vspace{16pt}
    \subfloat[\centering Instance 2: In-sample stability.] 
    {\label{fig:4p-iss} \includegraphics[width=8cm]{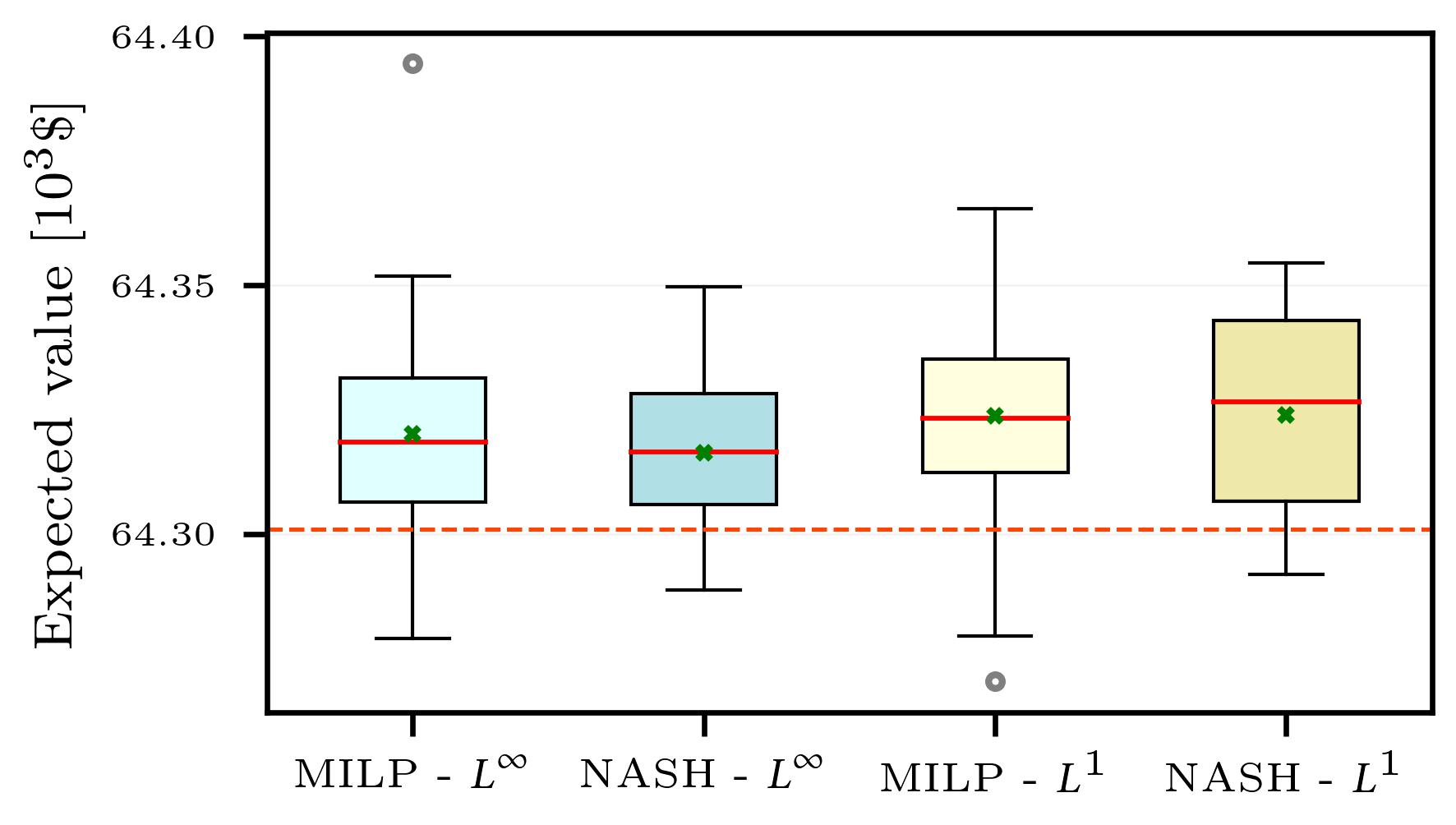} }%
    \subfloat[\centering Instance 2: Out-of-sample stability.] 
    {\label{fig:4p-oss} \includegraphics[width=8cm]{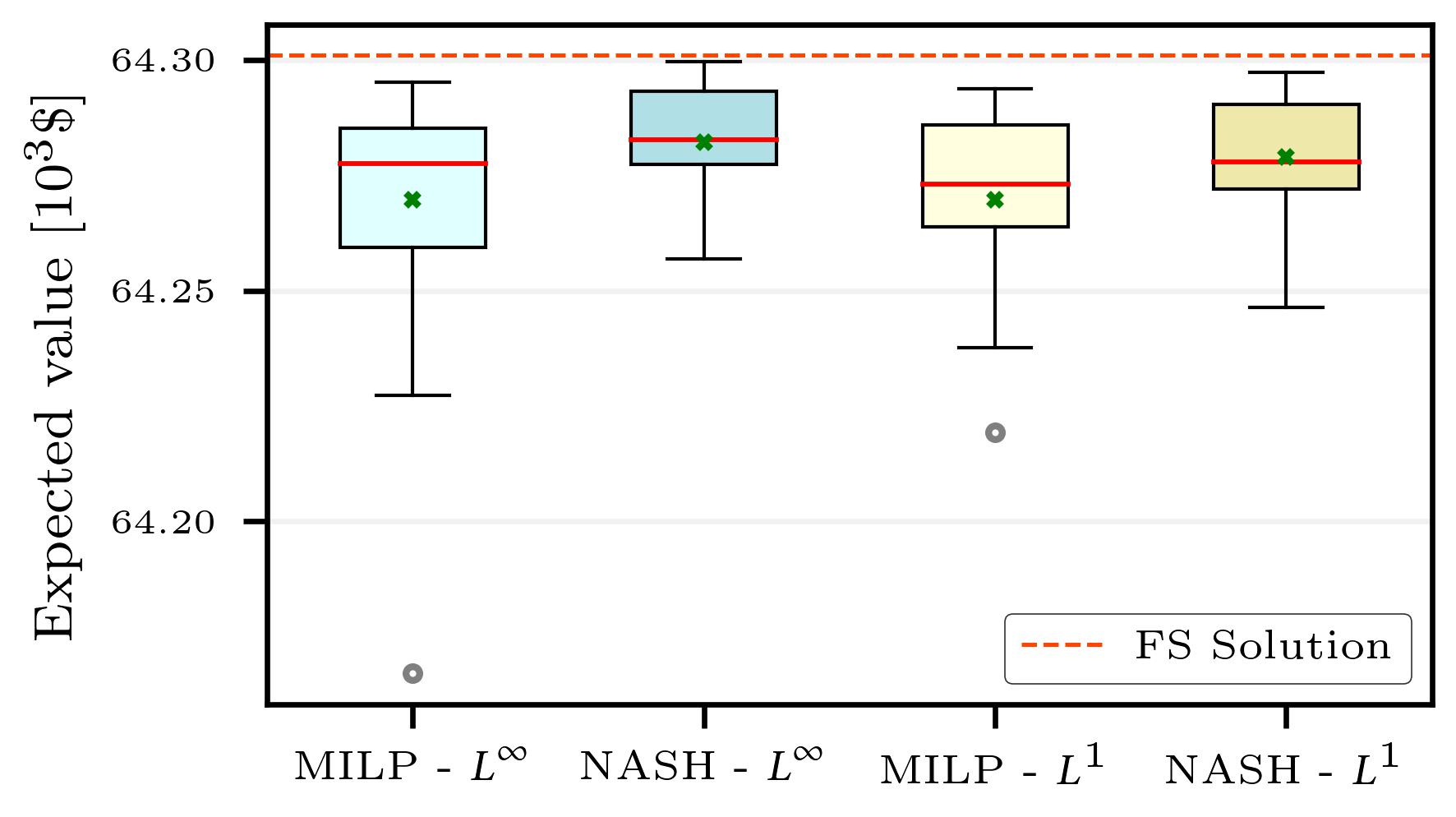} }%
    \\
    \vspace{16pt}
    \subfloat[\centering Instance 3: In-sample stability.] 
    {\label{fig:8p-iss} \includegraphics[width=8cm]{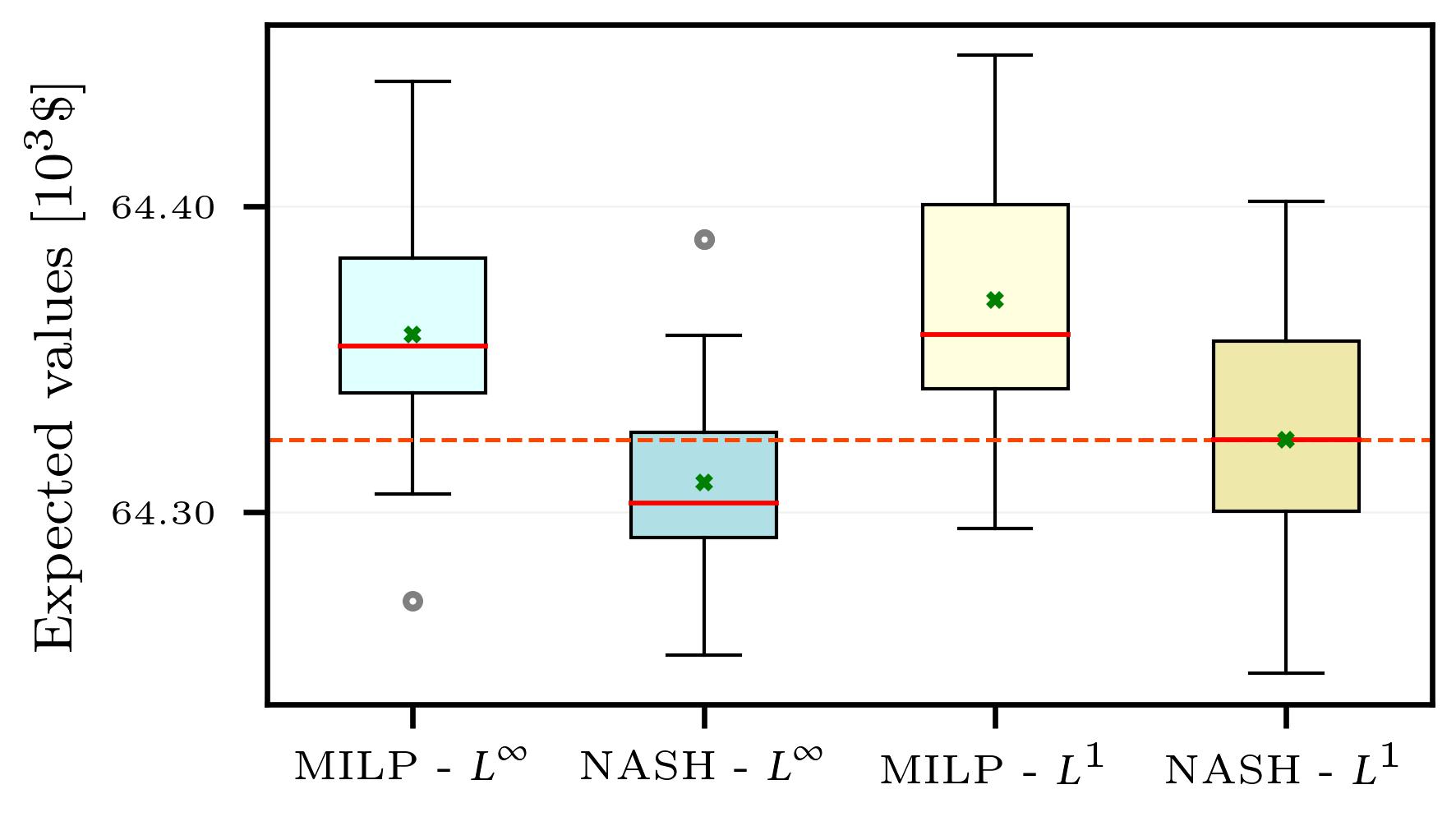} }%
    \subfloat[\centering Instance 3: Out-of-sample stability.] 
    {\label{fig:8p-oss} \includegraphics[width=8cm]{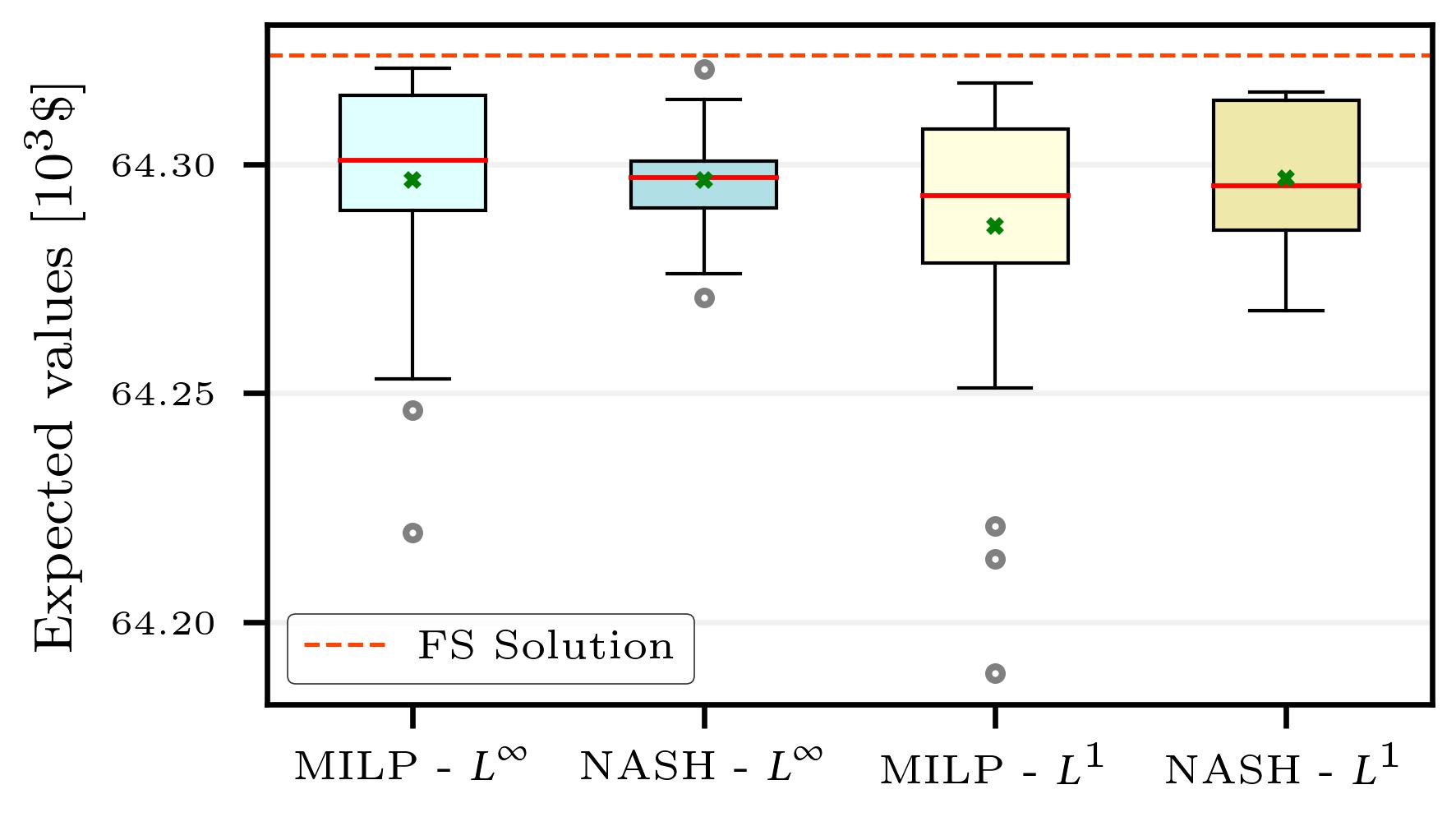} }%
    \caption{Stability tests for all instances of the case study.}%
    \label{fig:stability}%
\end{figure}

It is noted that this work differs from the existing approaches in the literature as the stability of the optimisation-based SG models is evaluated over the 37 sets of errors' weights. In general this approach is followed as the work by \citet{Bounitsis2022} have indicated variability of the results of different selections of errors' weights on the DMP models. Results on in-sample and out-of-sample stability for all instances are presented in Fig. \ref{fig:stability}. Moreover, it is mentioned that the results on Table \ref{table:bias} are representative of the bias and can be connected to the out-of-sample stability.

On the one hand, the results regarding the in-sample stability (Figs. \ref{fig:2p-iss}, \ref{fig:4p-iss}, \& \ref{fig:8p-iss}) demonstrate that both MILP and NASH models achieve similar levels of stability, while the $L^{\infty}$-norm-based models lead to slightly tighter interquartile ranges. In particular, for Instance 3 of 8 uncertain parameters the mean of the expected values achieved by the NASH models is close to the FS solution. As 40 scenario are generated to represent the uncertainty of 8 uncertain parameters in Instance 3, the in-sample results indicate that the scenario sets obtained from NASH models can consistently lead to more accurate expected values when used to solve the scenario-based TSSP.

On the other hand, the results on out-of-sample stability (Figs. \ref{fig:2p-oss}, \ref{fig:4p-oss}, \& \ref{fig:8p-oss}) are more appropriate to capture the accuracy of the stochastic solutions of the scenario-based TSSP on the full-space TSSP (employing the reference tree $R$). In these tests the results more vividly showcase the positive impact of the game theoretic extension for the performance of the DMP method. In particular, for all instances the NASH models display tighter interquartile ranges than the original MILP models. This indicates that the stochastic solutions obtained from the scenario-based TSSP impose a bias with lower variance on the full-space problem. Moreover, the mean values of the sets are consistently slightly closer to the full-space solution, indicating that apart from an increased stability of the NASH models the quality of the stochastic solution is improved using the NASH DMP models.

Ultimately, both MILP and NASH models display a very good stability for the instances of the investigated case study. The results regarding the MILP model can complement the work by \citet{Bounitsis2022} and show than not only enhanced quality of stochastic solutions is obtained but also this enhanced performance can be consistent over the selection of errors' weights. Especially the proposed game theoretic extension of the model is demonstrated to enhance the original MILP model. In other words, the proposed NASH DMP model may generate scenario sets of better quality with an increased level of confidence independently of the user-defined weights compared to the original MILP model.

\section{Concluding Remarks}
\label{sec:conclusions}

To summarise, this work studies optimisation-based scenario generation methods and it delves into the stability of these methods over user-defined parameters of the optimisation models. In particular, the Distribution and Moments Matching Problem (DMP) is studied as: (i) NLP and MILP optimisation models are already proposed for DMP \citep{Calfa2014,Bounitsis2022}, (ii) the objective function consists of several terms and user-defined parameters are indicated to induce variability in the solutions \citep{Bounitsis2022}. Thus, systematic measures assessing the stability and the quality of the stochastic solutions using SG methods are employed in this work over different sets of errors' weights. In comparison to the existing model a Nash game theoretic extension of the DMP MILP model (NASH) is proposed. This model considers the terms of the objective function as players of a game and is formulated as MILP following a separable programming approach. The results indicate that NLP may provide scenario sets of low quality because of high statistical errors and persisting under-specification issues. However, MILP and NASH models achieve acceptable statistical errors which generally lead to good stability and quality of stochastic solutions. Thus the results validate the stability of the DMP MILP models over the selection of user-defined weights. Moreover, NASH models further outperforms the MILP models regarding the stability and the quality of the stochastic solutions over various sets of user-defined errors' weights. Ultimately, the proposed NASH models is demonstrated by the analysis to constitute a more reliable approach towards scenario generation independently of the historical data at hand and the user-defined parameters.

The increased stability of the SG models and the avoidance of under-specification issues are crucial towards the incorporation of the single scenario reduction models to a multi-stage setting. In the future work within our group the goal is to exploit the high fidelity stable DMP MILP model in order to propose frameworks for multi-stage scenario generation using scenario lattice approach \citep{Bounitsis2023IN}. Finally, proposed frameworks are aimed to be used for the solution of real-world problems of the PSE domain under uncertainty towards explainable stochastic solutions \citep{Rathi2023}. 

\section*{CRediT authorship contribution statement }

\textbf{Georgios L. Bounitsis:} Conceptualisation, Methodology, Investigation, Formal analysis, Data Curation, Writing - Original Draft, Visualisation. \textbf{Lazaros G. Papageorgiou:} Conceptualisation, Review \& Editing, Supervision. \textbf{Vassilis M. Charitopoulos:} Conceptualisation, Methodology, Writing - Review \& Editing, Investigation, Supervision, Funding acquisition.

\section*{Declaration of competing interest }

The authors declare that they have no competing  interests.

\section*{Acknowledgements} 

Financial support from the UK EPSRC under projects EP/V034723/1, EP/V051008/1 \& EP/V050168/1 is gratefully acknowledged.

\section*{Nomenclature}
\label{sec:nomenclature}

\subsubsection*{Indices}
\begin{longtable*}{p{1.7cm} p{14cm}}
$i$,$i^\prime$ & uncertain parameters \\ 
$g$,$g^\prime$ & prespecified grid points \\ 
$k$,$k^\prime$ & final scenarios/clusters  \\ 
$m$,$m^\prime$ & statistical moments \\ 
$n$,$n^\prime$ & original scenarios/data points of the uncertain set \\ 
$t$,$t^\prime$ & players of the game/terms of objective function
\end{longtable*}

\subsubsection*{Sets}
\begin{longtable*}{p{1.7cm} p{14cm}}
$I$ & set of uncertain parameters \\ 
$G$ & set of prespecified grid points \\ 
$K$ & set of final scenarios/clusters\\ 
$M$ & set of statistical moments \\ 
$N$ & set of original scenarios/data points of the uncertain set \\
$T$ & set of players of the game \\
$CL_{kn}$ & subset of original scenarios $n$ which belong to cluster $k$ \\ 
\end{longtable*}

\subsubsection*{Parameters}
\begin{longtable*}{p{1.7cm} p{14cm}}
${C}_{ii^\prime}$ & covariance between uncertain parameters $i,i^\prime$\\ 
${D}_{im}$ & value of $m^{th}$ moment of uncertain parameter $i$ \\ 
${ECDF}_{in}$ & empirical cumulative probability of data point $n$ of uncertain parameter $i$\\ 
$P^{max}$ & maximum allowable probability of occurrence of a scenario \\
$P^{min}$ & minimum allowable probability of occurrence of a scenario \\ 
${X}_{in}$ & value of uncertain parameter $i$ in original scenario $n$ \\ 
$W^{COV}_{ii^\prime}$ & weight of the error regarding the covariance \\ 
$W^{ECDF}_{i}$ & weight of the error regarding the ECDF of uncertain parameter $i$ \\ 
$W^{SM}_{im}$ & weight of the error regarding the $m^{th}$  moment of one uncertain parameter $i$\\ 
$\overline{W}^{COV}_{ii^\prime}$ & arbitrarily chosen parameter for the calculation of $W^{COV}_{ii^\prime}$ \\ 
$\overline{W}^{ECDF}_{i}$ & arbitrarily chosen parameter for the calculation of $W^{ECDF}_{i}$ \\
$\overline{W}^{SM}_{im}$ & arbitrarily chosen parameter for the calculation of $W^{SM}_{im}$ \\ 
$\pi_t^{MAX}$ & maximum allowable error regarding player/term $t$  \\
$\overline{\pi}_{tg}$ & error regarding player/term $t$ at grid point $g$ \\ 
 
\end{longtable*}

\subsubsection*{Binary Variables}
\begin{longtable*}{p{1.7cm} p{14cm}}
$y_{kn}$ & is 1, if original scenarios $n$ is selected as final scenario $k$; 0, otherwise \\ 
$\lambda_{tg}$ & SOS2 variable for the approximation of the Nash product \\ 
\end{longtable*}

\subsubsection*{Continuous Variables}
\begin{longtable*}{p{1.7cm} p{14cm}}
${c}^{+}_{ii^\prime}$ & positive deviation on covariance between uncertain parameters $i,i^\prime$\\ 
${c}^{-}_{ii^\prime}$ & negative deviation on covariance between uncertain parameters $i,i^\prime$\\
${d}^{+}_{im}$ & positive deviation on value of $m^{th}$ moment of uncertain parameter $i$ \\ 
${d}^{-}_{im}$ & negative deviation on value of $m^{th}$ moment of uncertain parameter $i$ \\ 
${e}_{i}$ & maximum deviation on ECDF of uncertain parameter $i$ \\ 
$p_{kn}$ & probability of occurrence of original scenario $n$ when selected as final scenario $k$ \\ 
$\phi_{in}$ & deviation on ECDF of data point $n$ of uncertain parameter $i$ \\ 
$\widehat{\phi}_{ik}$ & deviation on ECDF of selected scenario $k$ of uncertain parameter $i$ \\ 
$\pi_t$ & error regarding player/term $t$  \\
\end{longtable*}

\section*{Appendices} 
\appendix

\section{Supplementary data}
\label{sec:appendix}
Mathematical notation regarding the problem of the case study in Section \ref{sec:casestudy} is reported as follows:

\subsubsection*{Indices}
\begin{longtable*}{p{2.2cm} p{13.5cm}}
$h$ & materials \\ 
$j$ & processes \\ 
$l$ & products
\end{longtable*}

\subsubsection*{Parameters}
\begin{longtable*}{p{2.2cm} p{13.5cm}}
$D_l$ & demand for product $l$ \\ 
${DC}_j$ & cost for process $j$ \\ 
${FC}_j$& fixed cost of process $j$ \\ 
${maxQ}_j$ & maximum volume capacity of process $j$\\
${maxRM}_h$ & maximum availability of raw material $h$ \\
${MI}_j$ & mass flow to volume relationship constant for process $j$ \\
${OC}_j$ & operating cost of process $j$ \\
${PC}_j$ & yield constant for process $j$ \\
$\alpha_h$ & cost of raw material $h$ \\
$\beta_l$ & price of product $l$
\end{longtable*}

\subsubsection*{Variables}
\begin{longtable*}{p{2.2cm} p{13.5cm}}
${IS}_j$ & mass flow in input stream to process $j$ \\ 
${OS}_j$ & mass flow in output stream to process $j$ \\ 
$P_l$& mass flow of product $l$ \\ 
$Q_j$ & capacity of process $j$ \\
${RM}_h$ & mass flow of raw material $h$ \\
$f_1,\ldots,f_5$ & mass flow of intermediate streams \\
$y_j$ & is 1, if process $j$ is selected; 0, otherwise
\end{longtable*}

In Table \ref{table:app:cs:data} are provided all the necessary data for the solution of the problem of the case study (Section \ref{sec:casestudy}). Values marked with asterisk symbol ($^{\ast}$) are used considered as uncertain parameters in some instances.

\begin{table}[H]
\renewcommand{\arraystretch}{1.5}
\renewcommand{\thetable}{A.1}
\fontsize{11pt}{11pt} \selectfont
\begin{center}
\caption{ Values of deterministic parameters of case study. }
\label{table:app:cs:data}
\begin{tabular}{l c c c c c c c c c c c } 

\hline
$\mathbf{h,j,l}$ & \textbf{1} & \textbf{2} & \textbf{3} & \textbf{4} & \textbf{5} & \textbf{6} & \textbf{7} & \textbf{8} & \textbf{9} & \textbf{10} & \textbf{11}   
\\	
\hline
$\mathbf{{DC}_j}$ & 2500 &	2500	&2500&	2500&	2500&	2500	&2500&	2500&	2500	&2500&	2500 \\
$\mathbf{{FC}_j}$ & 4000&	2500&	3500&	3000&	4500&	2500&	3000&	2200&	2800&	2700&	2500 \\
$\mathbf{{maxQ}_j}$ & 3.0 &	3.0	&3.0&	3.0&	3.0&	3.0&	3.0&	3.0&	3.0	&3.0&	3.0 \\
$\mathbf{{MI}_j}$ & 18&	20&	15&	20&	20&	21&	15&	15&	25&	15&	20 \\
$\mathbf{{OC}_j}$ & 400&	400&	400&	400&	400&	400&	400&	400&	400&	400&	400\\
$\mathbf{{PC}_j}$ & 13&	15&	17&	14$^{\ast}$&	10&	15&	16$^{\ast}$ &	11$^{\ast}$ &	13&	15&	17$^{\ast}$\\
$\mathbf{D_l}$ & 30.80&	29.60$^{\ast}$ &	30.05$^{\ast}$ &	29.50$^{\ast}$ &	30.00$^{\ast}$ & \multicolumn{6}{c}{ } \\
$\mathbf{{maxRM}_h}$ & 34.80&	35.65&	33.65&	35.50&	35.00 & \multicolumn{6}{c}{ } \\
$\mathbf{{\boldsymbol{\alpha}}_h}$ & 200&	320&	230&	250&	300 & \multicolumn{6}{c}{ } \\
$\mathbf{{\boldsymbol{\beta}}_l}$ & 600&	650&	500&	400&	700 & \multicolumn{6}{c}{ } \\
\hline

\end{tabular}
\end{center}
\end{table}

\bibliographystyle{elsarticle-harv} 
\bibliography{draft-dmpnash-abbr}





\end{document}